\documentclass[draft,reqno,12pt]{amsart}

\newtheorem{theorem}{Theorem}[section]
\newtheorem{proposition}[theorem]{Proposition}
\newtheorem{lemma}[theorem]{Lemma}
\newtheorem{cor}[theorem]{Corollary}

\theoremstyle{definition}

\newtheorem{exam}[theorem]{Example}
\newtheorem{problem}[theorem]{Problem}
\newtheorem*{theo1}{Theorem \ref{soperator}}

\newcommand{\Gin}{\ensuremath{\mathrm{Gin}}}
\newcommand{\gin}{\ensuremath{\mathrm{gin}}}

\newcommand{\init}{\ensuremath{\mathrm{in}}}
\newcommand{\GL}{\ensuremath{GL_n (K)}}

\newcommand{\dele}[1]{\ensuremath{\Delta^e({#1})}}

\newcommand{\rlex}{{<_{\mathrm{rev}}}}
\newcommand{\aru}[1]{\ensuremath{R_{[{#1}]}}}
\newcommand{\setmono}{\ensuremath{ M_{[\infty]}}}

\begin{document}
\title{Generic initial ideals and squeezed spheres}
\author{Satoshi Murai}
\address{
Department of Pure and Applied Mathematics\\
Graduate School of Information Science and Technology\\
Osaka University\\
Toyonaka, Osaka, 560-0043, Japan\\
}
\email{s-murai@ist.osaka-u.ac.jp
}
\thanks{The author is supported by JSPS Research Fellowships for Young Scientists}
\maketitle

\begin{abstract}
In 1988 Kalai construct a large class of simplicial spheres,
called squeezed spheres,
and in 1991 presented a conjecture about generic initial ideals of 
Stanley--Reisner ideals of squeezed spheres.
In the present paper this conjecture will be proved.
In order to prove Kalai's conjecture,
based on the fact that every squeezed $(d-1)$-sphere is the boundary of a 
certain $d$-ball,
called a squeezed $d$-ball,
generic initial ideals of Stanley--Reisner ideals of squeezed balls will be 
determined.
In addition, generic initial ideals of exterior face ideals of squeezed 
balls are determined.
On the other hand, we study the squeezing operation,
which assigns to each Gorenstein* complex $\Gamma$
having the weak Lefschetz property a squeezed sphere $\mathrm{Sq}(\Gamma)$,
and show that this operation increases graded Betti numbers.
\end{abstract}

\section*{Introduction}

Let $K$ be an infinite field,
$ \aru n=K[x_1,x_2,\dots,x_n]$
the polynomial ring over a field $K$ with each $\deg(x_i)=1$.
For a graded ideal $I \subset \aru n$,
let $\gin(I)$ be the generic initial ideal of $I$ with respect to the degree reverse lexicographic order 
induced by $x_1>x_2>\cdots >x_n$.
Let $u=x_{i_1} x_{i_2}\cdots x_{i_d} \in \aru n $ 
and $v = x_{j_1} x_{j_2} \cdots x_{j_d} \in \aru n$
be monomials of degree $d$ with $i_1 \leq i_2 \leq \cdots \leq i_d$ and with $j_1 \leq j_2 \leq \cdots \leq j_d$.
We write $u \prec v$ if $ i_k \leq j_k$ for all $1\leq k \leq d$.
A monomial ideal $I\subset \aru n$ is called {\em strongly stable}
if $v \in I$ and $u \prec v$ imply $u \in I$. 
Generic initial ideals are strongly stable if the base field is of characteristic $0$.

Applying the theory of generic initial ideals to combinatorics
by considering generic initial ideals of Stanley--Reisner ideals or exterior face ideals
is known as algebraic shifting.
Kalai \cite{K} proposed a lot of problems about algebraic shifting.
In the present paper,
we will prove a problem in \cite{K} about generic initial ideals of Stanley--Reisner ideals of squeezed spheres.

Squeezed spheres were introduced by Kalai \cite{K3} by extending the construction of Billera--Lee polytopes.
For a simplicial complex $\Gamma$ on the vertex set $[n]=\{1,2,\dots,n\}$,
let $I_\Gamma$ be the Stanley--Reisner ideal of $\Gamma$.
Fix integers $n>d>0$ and $m \geq 0$.
A set $U \subset \aru m$ of monomials is called
{\em a shifted order ideal of monomials} if $U$ satisfies
\begin{itemize}
\item[(i)]
$\{1,x_1,x_2,\dots,x_m\}\subset U$;
\item[(ii)]
if $u \in U$ and $v \in \aru m$ divides $u$, then $v \in U$;
\item[(iii)]
if $u\in U$ and $u \prec v$, then $v \in U$.
\end{itemize}
(In general, (i) is not assumed.
The reason why we assume it will be explained in \S 2.)
If  $U \subset \aru {n-d-1}$ is a shifted order ideal of monomials 
of degree at most $\lfloor \frac{d+1} 2 \rfloor$,
then we can construct a shellable $d$-ball $B_d(U)$ on $[n]$
by considering a certain subcollection of the collection of facets of the boundary complex of the cyclic $d$-polytope with $n$ vertices
associate with $U$.
The squeezed sphere $S_d(U)$ is the boundary of the squeezed $d$-ball $B_d(U)$.
The $g$-vectors of squeezed spheres are given by 
$g_i(S_d(U))=|\{u \in U:\deg(u)=i\}|$ for $ 0 \leq i \leq \lfloor \frac d 2 \rfloor$,
where $|V|$ denotes the cardinality of a finite set $V$.

On the other hand,
if $\Gamma$ is a $(d-1)$-dimensional Gorenstein* complex on $[n]$ with the weak Lefschetz property
and if the base field $K$ is of characteristic $0$,
then the set of monomials
$$U(\Gamma)=\{ u \in \aru {n-d-1}: u\not \in \gin(I_\Gamma) \mbox{ is a monomial} \}$$
is a shifted order ideal of monomials of degree at most $\lfloor \frac d 2 \rfloor$
with $g_i(\Gamma)=|\{ u\in U(\Gamma):\deg(u)=i\}|$.
Furthermore, if $\Gamma$ has the strong Lefschetz property,
then $U(\Gamma)$ determine $\gin(I_\Gamma)$.
(See \S 3.)

In \cite{K2} and \cite[Problem 24]{K},
Kalai conjectured that,
for any shifted order ideal $U \subset \aru {n-d-1}$ of monomials of degree at most $\lfloor \frac d 2 \rfloor$, one has
\begin{eqnarray}
U(S_d(U))=U. \label{intro}
\end{eqnarray}
In the present paper, the above conjecture will be proved affirmatively
(Theorem \ref{conj}).

To solve the above problem, the concept of stable operators is required.
Let $\aru \infty=K[x_1,x_2,x_3,\dots]$ be the polynomial ring in infinitely many variables
and $\setmono$ the set of monomials in $\aru \infty$.
Let $\sigma : \setmono \to \setmono$ be a map,
$I\subset \aru \infty$ a finitely generated monomial ideal
and $G(I)$ the set of minimal monomial generators of $I$.
Write $\sigma(I)$ for the ideal generated by $\{\sigma(u):u\in G(I)\}$.
Let $\beta_{ij}(I)$ be the graded Betti numbers of the ideal $I \cap \aru m$ and $\gin(I)=\gin(I \cap \aru m) \aru \infty$
for an integer $m$ with $G(I) \subset \aru m$.
A map $\sigma: \setmono \to \setmono$ is called a {\em stable operator}
if $\sigma$ satisfies
\begin{itemize} 
\item[(i)]
if $I\subset \aru \infty$ is a finitely generated strongly stable ideal,
then $\beta_{ij}(I)=\beta_{ij}(\sigma (I))$ for all $i,j$;
\item[(ii)]
if $J\subset I$ are finitely generated strongly stable ideals of $\aru \infty$,
then $\sigma (J) \subset \sigma (I)$.
\end{itemize}
The first result is the following.
(Similar results can be found in \cite{BNT} and \cite{BCR}.)

\begin{theo1} 
Let $\sigma:\setmono \to \setmono$ be a stable operator.
If $I\subset \aru \infty$ is a finitely generated strongly stable ideal,
then $\gin(\sigma(I))=I$.
\end{theo1}



By using Theorem \ref{soperator},
we determine generic initial ideals of Stanley--Reisner ideals of squeezed balls 
(Proposition \ref{ginball}).
A squeezed sphere $S_d(U)$ is called an \textit{S-squeezed $(d-1)$-sphere} if $U\subset \aru {n-d-1}$ is a shifted order ideal of monomials
of degree at most $\lfloor \frac d 2 \rfloor$.
If $S_d(U)$ is an S-squeezed sphere,
then $B_d(U)$ and $S_d(U)$ have the same $\lfloor \frac {d+1} 2 \rfloor$-skeleton.
By using this fact together with the forms of $\gin(I_{B_d(U)})$,
we will show the equality (\ref{intro}).
In particular, the equality (\ref{intro}) immediately implies that
every S-squeezed sphere has the weak Lefschetz property.

This paper is organized as follows.
In \S 1, we will study stable operators.
In \S 2 and \S 3, we recall some basic facts about squeezed spheres and Lefschetz properties.
In \S 4, we will prove Kalai's conjecture.
We also study some properties of S-squeezed spheres in \S 5, \S 6 and \S 7.

In \S 5, we study the relation between squeezing and graded Betti numbers.
Assume that the base field is of characteristic $0$.
Let $\Gamma$ be a $(d-1)$-dimensional Gorenstein* complex with the weak Lefschetz property.
Then $U(\Gamma)=\{ u \in \aru {n-d-1}: u\not \in \gin(I_\Gamma) \mbox{ is a monomial} \}$
is a shifted order ideal of monomials of degree at most $\lfloor \frac d 2 \rfloor$.
Define $\mathrm{Sq}(\Gamma) = S_d(U(\Gamma))$.
Then $\Gamma $ and $\mathrm{Sq}(\Gamma)$ have the same $f$-vector.
Also, by virtue of the equality (\ref{intro}),
we have $\mathrm{Sq}(\mathrm{Sq}(\Gamma)) = \mathrm{Sq}(\Gamma)$.
Although generic initial ideals do not preserve the Gorenstein property,
we can define the operation $\Gamma \to \mathrm{Sq}(\Gamma)$ which assigns to each Gorenstein* complex $\Gamma$ 
having the weak Lefschetz property an S-squeezed sphere $\mathrm{Sq}(\Gamma).$
This operation $\Gamma \to \mathrm{Sq}(\Gamma)$ is called \textit{squeezing}.

First, we will show that the graded Betti numbers of the Stanley--Reisner ideal of each $\mathrm{Sq}(\Gamma)$ 
are easily computed by using Eliahou--Kervaire formula together with $U(\Gamma)$ (Theorem \ref{sqbetti}).
Second,
we will show that squeezing increases graded Betti numbers
(Theorem \ref{sqinc}).


In \S 6,
we consider S-squeezed 4-spheres.
Since Pfeifle proved that squeezed $3$-spheres are polytopal,
we can easily show that S-squeezed $4$-spheres are polytopal.
This fact yields a complete characterization of generic initial ideals of Stanley--Reisner ideals
of the boundary complexes of simplicial $d$-polytopes for $d \leq 5$,
when the base field is of characteristic $0$.


In \S 7, we consider generic initial ideals in the exterior algebra.
We will determine generic initial ideals of exterior face ideals of squeezed balls
by using the technique of squarefree version of stable operators.

\section{generic initial ideals and stable operators}
Let $K$ be an infinite field,
$\aru n =K[x_1,x_2,\dots,x_n]$
the polynomial ring in $n$ variables over a field $K$ with each $\deg(x_i)=1$
and $M_{[n]}$ the set of monomials in $\aru n$.
Let $\aru \infty=K[x_1,x_2,x_3,\dots]$ be the polynomial ring in infinitely many variables
and $\setmono$ the set of monomials in $\aru \infty$.
For a graded ideal $I \subset \aru n$ and for an integer $d \geq 0$,
let $I_d$ denote the $d$-th homogeneous component of $I$.

Fix a term order $<$ on $\aru n$.
For any polynomial $f= \sum_{u \in M_{[n]}} \alpha_u u \in \aru n$
with each $\alpha _u \in K$,
the monomial $\init_{<}(f)= \max_{<} \{ u: \alpha_u \ne 0\}$
is called the {\em initial monomial of $f$}.
The \textit{initial ideal} $ \init_< (I)$ of an ideal $I\subset \aru n$
is the monomial ideal generated by
the initial monomials of all polynomials in $I$.


Let $\GL$ be the general linear group with coefficients in $K$.
For any $\varphi=(a_{ij}) \in \GL$ and for any polynomial $f\in \aru n$,
define
$$ \varphi (f(x_1,x_2,\dots,x_n))=f(\sum_{i=1}^n a_{i1}x_i,\sum_{i=1}^n a_{i2}x_i,\dots,\sum_{i=1}^n a_{in}x_i).$$
For a graded ideal $I$,
we let $\varphi(I)=\{\varphi(f):f\in I\}$.

The fundamental theorem of generic initial ideals is the following.
\begin{theorem}[Galligo, Bayer and Stillman] \label{gin}
Fix a term order $<$ satisfying $x_1 > x_2 > \dots > x_n$.
For each graded ideal $I \subset \aru n$,
there is a nonempty Zariski open subset $U\subset \GL$ such that
$\init_{<} ( \varphi (I))$ is constant for all $\varphi \in U$.
Furthermore, if $K$ is a field of characteristic $0$, then $\init_{<} (\varphi (I))$ with $\varphi \in U$ is strongly stable.
\end{theorem}

This monomial ideal $\init_{<}(\varphi (I))$ with $\varphi \in U$ is called the
{\em generic initial ideal of $I$ with respect to the term order $<$},
and will be denoted $\gin_{<} (I)$.
Let $\rlex$ be the degree reverse lexicographic order induced by $x_1>x_2>x_3>\cdots$.
In other words,
for monomials $u=x_{i_1} x_{i_2}\cdots x_{i_k} \in \aru \infty$ 
and $v = x_{j_1} x_{j_2} \cdots x_{j_l} \in \aru \infty$
with $i_1 \leq i_2 \leq \cdots \leq i_k$ and with $j_1 \leq j_2 \leq \cdots \leq j_l$,
one has $ u\rlex v$ if $\deg(u) < \deg (v)$ or $\deg(u)=\deg(v)$ and
for some $r$ one has $i_t =j_t$ for $t > r$, and $i_r>j_r$.
In the present paper,
we only consider generic initial ideals w.r.t.\ the degree reverse lexicographic order and
write $\gin (I)= \gin_{\rlex} (I).$
We recall some fundamental properties.

\begin{lemma}[{\cite[Lemma 3.3]{BNT}}] \label{bnt}
Let $I \subset \aru n$ be a graded ideal.
Then
\begin{itemize}
\item[(i)]
$I$ and $\gin (I)$ have the same Hilbert function.
In other words,
$\dim_K (I_d)=\dim_K (\gin (I)_d)$ for all $d \geq 0$;
\item[(ii)]
if $J \subset I$ are graded ideals of $\aru n$,
then $\gin(J) \subset \gin(I)$;
\item[(iii)]
$\gin (I \aru {n+1})= (\gin(I))\aru {n+1}$.
\end{itemize}
\end{lemma}

Let $M$ be a finitely generated graded $\aru n$-module.
The {\em graded Betti numbers} $\beta_{ij}=\beta_{ij}(M)$ of $M$,
where $i,j \geq 0$,
are the integers $\beta_{ij} (M) = \dim_K (\mathrm{Tor}_i(M,K)_j)$.
In other words,
$\beta_{ij}$ appears in the minimal graded free resolution
\[
0 \longrightarrow \bigoplus_{j} \aru n(-j)^{\beta_{hj}}
\longrightarrow
\cdots
\longrightarrow
\bigoplus_{j} \aru n (-j)^{\beta_{1j}}
\longrightarrow
\bigoplus_{j} \aru n (-j)^{\beta_{0j}}
\longrightarrow M
\longrightarrow 0
\]
of $M$ over $\aru n$.
The {\em projective dimension of $M$} is the integer
$$\mathrm{proj\ dim} (M)= \max\{ i: \beta_{ij}(M) \ne 0 \mbox{ for some }j\geq 0 \}.$$

\begin{lemma}[{\cite[Corollary 19.11]{E}}] \label{extremal}
Let $I\subset \aru n$ be a graded ideal.
Then
$$\mathrm{proj \ dim}(I)=\mathrm{proj \ dim} (\gin(I)).$$
\end{lemma}

For any monomial $u \in \setmono$,
write ${m}(u)=\max\{i: x_i \mbox{ divides } u\}$.
Recall that every generic initial ideal $\gin(I)$ of a graded ideal $I \subset \aru n$ is Borel-fixed,
that is,
one has $\varphi(\gin(I)) =\gin (I)$ for any upper triangular invertible matrix $\varphi \in \GL$
(see \cite[Theorem 15.20]{E}).  
It follows from \cite[Corollary 15.25]{E} together with the Auslander--Buchsbaum formula (\cite[Theorem 19.9]{E})
that the projective dimension of any Borel-fixed monomial ideal $J \subset \aru n$ is
$$\mathrm{ proj \ dim}(J)=\mathrm{proj \ dim}(\aru n/J)-1= \max\{ m (u): u \in G(J)\}-1,$$ 
where $G(J)$ is the set of minimal monomial generators of $J$.
Thus the next lemma immediately follows from Lemma \ref{extremal} together with the above equality.
 
\begin{lemma} \label{idealmax}
Let $I \subset \aru n$ be a graded ideal.
Then 
$$\mathrm{proj \ dim}(I)= \max\{ m(u): u\in G(\gin(I))\}-1.$$
\end{lemma}

Let $I \subset \aru \infty$ be a finitely generated monomial ideal and
$G(I)$ the set of minimal monomial generators of $I$.
Write $\max(I) = \max \{ m(u): u \in G(I)\}$.
We let $\gin(I) = \gin(I\cap \aru {\max(I)})\aru \infty$
and let $\beta_{ij}(I)$ (resp. $\mathrm{proj \ dim}(I)$)
be the graded Betti numbers (resp. the projective dimension)
 of $I\cap \aru {\max(I)}$ over $\aru {\max(I)}$.
Note that the graded Betti numbers of $I\cap \aru {k}$ over $\aru k$ are constant for all $k \geq \max (I)$.
Also, Lemma \ref{bnt} (iii) guarantees
$\gin(I\cap \aru k)\aru \infty=
\gin( I \cap \aru {\max(I)}) \aru \infty$ for all $k \geq \max(I)$.
We say that a finitely generated monomial ideal $I \subset \aru \infty$ is strongly stable
if $I \cap \aru {\max(I)}$ is strongly stable.

An important fact on generic initial ideals
is that graded Betti numbers of strongly stable ideals are easily computed
by the Eliahou--Kervaire formula.
We recall the Eliahou--Kervaire formula.

\begin{lemma}[{\cite[Corollary 3.4]{H}}] \label{e-k}
Let $I \subset \aru n$ be a strongly stable ideal.
Then
\begin{itemize}
\item[(i)]
$\beta_{ii+j} (I)=\sum_{u \in G(I),\ \deg(u)=j} { {m}(u) -1 \choose i}$;
\item[(ii)]
$\mathrm{proj\ dim}(I)= \max(I)-1$.
\end{itemize}
\end{lemma}

Let $\sigma:\setmono \to \setmono$ be a map,
$I \subset \aru \infty$ a finitely generated monomial ideal and
$G(I)= \{ u_1,u_2,\dots,u_m\}$ the set of minimal monomial generators of $I$.
We write $\sigma (I)\subset \aru \infty$ for the monomial ideal
generated by $\{ \sigma(u_1),\sigma(u_2),\dots,\sigma(u_m) \}$. 

We say that
a map $\sigma: \setmono \to \setmono$ is a {\em stable operator}
if $\sigma$ satisfies
\begin{itemize} 
\item[(i)]
if $I\subset \aru \infty$ is a finitely generated strongly stable ideal,
then $\beta_{ij}(I)=\beta_{ij}(\sigma (I))$ for all $i,j$;
\item[(ii)]
if $J\subset I$ are finitely generated strongly stable ideals of $\aru \infty$,
then $\sigma (J) \subset \sigma (I)$.
\end{itemize}

\begin{theorem} \label{soperator}
Let $\sigma:\setmono \to \setmono$ be a stable operator.
If $I\subset \aru \infty$ is a finitely generated strongly stable ideal,
then $\gin(\sigma(I))=I$.
\end{theorem}

\begin{proof}
Let $m=\max(I)$.
Since $I$ is strongly stable,
Lemma \ref{e-k} says 
$\mathrm{proj \ dim}(I)=\max(I)-1$.
Also, since $I$ and $\sigma(I)$ have the same graded Betti numbers,
Lemma \ref{idealmax} says
$$\max(I)-1=\mbox{proj dim}(I)=\mathrm{proj \ dim}(\sigma(I))=\max(\gin(\sigma(I)))-1.$$
Then we have $\max(I)=\max(\gin(\sigma(I)))$.
Thus what we must prove is $ \gin(\sigma(I)) \cap \aru m= I\cap \aru m$.

We claim $I\cap \aru m$ and $\gin(\sigma(I)) \cap \aru m $ have the same Hilbert function.
Let $n = \max(\sigma(I))$.
We remark $m=\max(\gin (\sigma(I))) \leq n$.
Since 
$I$ and $\sigma(I)$ have the same graded Betti numbers,
$I\cap \aru n$, $\sigma(I) \cap \aru n$ and $\gin (\sigma (I)) \cap \aru n$ have the same Hilbert function.
Since $\max(I)=\max(\gin(\sigma(I)))=m$,
it follows that
$I\cap \aru m$ and $\gin(\sigma(I)) \cap \aru m$
have the same Hilbert function.

Now, we will show $ \gin(\sigma(I)) \cap \aru m= I\cap \aru m$ by using induction on $m$.
If $m=1$, then $G(I)$ is of the form $G(I)= \{x_1^k\}$,
where $k > 0$ is a positive integer.
Since $\max(\gin(\sigma(I)))=\max(I)=1$ and since 
$I\cap K[x_1]$ and $\gin(\sigma(I)) \cap K[x_1]$ have the  same Hilbert function,
we have $G(\gin(\sigma(I)))=\{x_1^k\}$.

Assume $m>1$.
Fix an integer $d\geq 0$.
Let  $I_{(d)}\subset \aru \infty$ be the ideal generated by all monomials $u \in I \cap \aru m$ of degree $d$.
Consider the ideal
$$J=(I_{(d)}:x_m^\infty)=\{f\in \aru \infty: ^\exists k\geq 0 \mbox{ such that } x_m^kf\in I_{(d)}\}.$$
Then $J$ is a finitely generated strongly stable ideal with $\max(J) <m$.

We claim 
\begin{eqnarray}
J_d \cap \aru m = (I_{(d)})_d \cap \aru m = I_d \cap \aru m. \label{e1}
\end{eqnarray}
Since $J \supset I_{(d)}$ and $(I_{(d)} )_d\cap \aru m= I_d \cap \aru m$ are obvious,
we will show $J_d \cap \aru m\subset (I_{(d)})_d\cap \aru m$.

Let
$ux_m^l \in J_d\cap \aru m$ be a monomial with $u \in \aru {m-1}$.
By the definition of $J=(I_{(d)}:x_m^\infty)$,
there is an integer $k \geq 0$ and a monomial $v x_m^{d-\deg(v)} \in G(I_{(d)})$
with $v \in \aru {m-1}$ such that $vx_m^{d-\deg(v)}$ divides $ux_m^{l+k}$.
This fact says $ux_m^l \prec v x_m^{d-\deg(v)}$.
Since $I_{(d)}$ is strongly stable,
we have $u x_m^l \in I_{(d)}$.
Thus we have $J_d \cap \aru m = (I_{(d)})_d \cap \aru m= I_d \cap \aru m$.

Since $I \supset I_{(d)}$ are strongly stable ideals,
Lemma \ref{bnt} together with the definition of stable operators says
\begin{eqnarray}
\gin(\sigma (I)) \supset \gin (\sigma (I_{(d)})). \label{hosi11}
\end{eqnarray}
Also, since $J \supset I_{(d)}$ are strongly stable ideals
and since $\max(J)<m$, the assumption of induction says
\begin{eqnarray}
J=\gin(\sigma (J)) \supset \gin (\sigma(I_{(d)})). \label{hosi12}
\end{eqnarray}
We already proved that if $I' \subset \aru \infty$ is a finitely generated strongly stable ideal
with $\max(I') \leq m$,
then $I' \cap \aru m$ and $\gin (\sigma (I')) \cap \aru m$ have the same Hilbert function.
This fact together with (\ref{e1}) says
$$\dim_K ( \gin( \sigma(J) )_d\cap {\aru m})=\dim_K(\gin (\sigma(I_{(d)}) )_d\cap {\aru m})=\dim_K (\gin (\sigma(I))_d\cap {\aru m}).$$
The above equality together with (\ref{e1}), (\ref{hosi11}) and (\ref{hosi12}) says
$$I_d \cap \aru m \stackrel{}{=}  J_d \cap \aru m\stackrel{}{=}
\gin(\sigma(J))_d \cap \aru m \stackrel{\star 1}{=} \gin(\sigma (I_{(d)}))_d \cap \aru m \stackrel{\star 2}{=} \gin(\sigma(I))_d \cap \aru m,$$
where ($\star 1$) follows from the inclusion (\ref{hosi12}) together with the fact that 
$\gin( \sigma(J) )_d\cap {\aru m}$ and $\gin( \sigma(I_{(d)}) )_d\cap {\aru m}$ are $K$-vector spaces with the same dimension
(and the equality ($\star2$) follows from the inclusion (\ref{hosi11}) by the same way as ($\star1$)).
Thus we have $I_d \cap \aru m = \gin(\sigma(I))_d \cap \aru m$ for all $d \geq 0$,
and therefore we have $I\cap \aru m = \gin(\sigma(I) )\cap \aru m$ as required.
\end{proof}

We will introduce an example of stable operators.
Let $a=(0,a_1,a_2,a_3,\dots)$ be a nondecreasing infinite sequence of integers.
Define the map $\alpha^a:\setmono \to \setmono$ by
\begin{eqnarray}
\alpha^a(x_{i_1}x_{i_2}\cdots x_{i_k})= x_{i_1} x_{i_2 + a_1} x_{i_3 + a_2} \cdots x_{i_k + a_{k-1}}, \label{teigi}
\end{eqnarray}
for any monomial $x_{i_1}x_{i_2}\cdots x_{i_k} \in \setmono$ with $i_1 \leq i_2 \leq \cdots \leq i_k$.
This map $\alpha^a$ is a generalization of the map studied in \cite{SO}.
We will show that the map $\alpha^a: \setmono \to \setmono $ is a stable operator.

Let $I \subset \aru \infty$ be a finitely generated monomial ideal,
$G(I)$ the set of minimal monomial generators of $I$ and $n = \max(I)$.
The ideal $I$ is said to have \textit{linear quotients} if for some order $u_1,u_2,\dots,u_m$ of the elements of $G(I)$
and for all $j=1,2,\dots,m$,
the colon ideals 
$$(\langle u_1,u_2,\dots,u_{j-1}\rangle :u_j)=\{ f \in \aru \infty : f u_j \in \langle u_1,u_2,\dots,u_{j-1} \rangle\}$$
are generated by a subset of $\{x_1,x_2,\dots,x_n\}$,
where $\langle u_1,u_2,\dots,u_{j-1} \rangle $ denotes the ideal generated by $\{ u_1,u_2,\dots,u_{j-1}\}$.
Define 
$$\mathrm{set}(u_j)=\{k \in [n]:x_{k} \in ( \langle u_1,u_2,\dots,u_{j-1} \rangle : u_j)\}\hspace{20pt} \mbox{ for } j=1,2,\dots,m.$$
If a finitely generated monomial ideal $I \subset \aru \infty$  has linear quotients,
then the graded Betti numbers of $I$ are given by the formula (\cite[Corollary 1.6]{HT})
\begin{eqnarray}
\beta_{ii+j}(I) = \sum _{u \in G(I),\ \deg(u)=j} { |\mathrm{set}(u)| \choose i}. \label{hosilq}
\end{eqnarray}

\begin{lemma} \label{include}
Let $a=(0,a_1,a_2,a_3,\dots)$ be a nondecreasing infinite sequence of integers and
$\alpha^a:\setmono \to \setmono$ the map defined in (\ref{teigi}).
Let $I\subset \aru \infty$ be a finitely generated strongly stable ideal.
If $u\in I$, then $\alpha^a(u) \in \alpha^a(I)$.
\end{lemma}

\begin{proof}
Let $u=x_{i_1}x_{i_2}\cdots x_{i_k} \in I$ with $i_1 \leq i_2 \leq \cdots \leq i_k$.
Since $u \in I$, there is $w \in G(I)$ such that $w$ divides $u$.
Since $I$ is strongly stable,
we may assume $w= x_{i_1}x_{i_2}\cdots x_{i_l}$ for some $l \leq k$.
Then $\alpha^a(w)=x_{i_1}x_{i_2 +a_1}\cdots x_{i_l+a_{l-1}} \in G(\alpha^a(I))$
divides $\alpha^a(u)=x_{i_1}x_{i_2 +a_1}\cdots x_{i_l+a_{l-1}}x_{i_{l+1}+a_l}\cdots x_{i_k+a_{k-1}}$.
Thus $\alpha^a(u) \in \alpha^a(I)$. 
\end{proof}


Let $u=x_1^{a_1}x_2^{a_2}\cdots x_k^{a_k}$ and $v=x_1^{b_1}x_2^{b_2}\cdots x_k^{b_k}$ be monomials with $m (u) \leq k$
and with $m(v) \leq k$.
The \textit{lexicographic order} $<_{\mathrm{lex}}$ of $\aru \infty$ is the total order on $\setmono$ defined by
$u <_{\mathrm{lex}} v$ if the leftmost nonzero entry of $(b_1-a_1,b_2-a_2,\dots,b_k-a_k)$ is positive.

\begin{lemma} \label{quotients}
With the same notation as in Lemma \ref{include}.
Let $G(I)=\{u_1,u_2,\dots,u_m\}$ be
the set of minimal monomial generators of $I$ with $u_1 >_{\mathrm{lex}} u_2 >_{\mathrm{lex}} \dots >_{\mathrm{lex}} u_m$.
Then $\alpha^a (I)\subset \aru \infty$ has linear quotients for the order 
$\alpha^a(u_1),\alpha^a(u_2),\dots,\alpha^a(u_m)$
with
$$\mathrm{set}(\alpha^a(x_{i_1}x_{i_2}\cdots x_{i_d} ))= \bigcup_{l=0}^{d-1} 
\{ k \in \mathbb{Z}: i_l + a_l \leq k < i_{l+1} +a_l \}$$
for any $x_{i_1}x_{i_2}\cdots x_{i_d} \in G(I)$ with $i_1 \leq i_2 \leq \cdots \leq i_d$,
where we let $i_0=1$ and $a_0=0$.
\end{lemma}

\begin{proof}
Set $u_j =x_{i_1}x_{i_2}\cdots x_{i_d} \in G(I) $ with $i_1 \leq i_2 \leq \cdots \leq i_d$
and
$$A(\alpha^a,u_j)=\bigcup_{l=0}^{d-1} 
\{ k \in \mathbb{Z}: i_l + a_l \leq k < i_{l+1} +a_l \}.$$

First, we will show
$\{x_k: k\in A(\alpha^a,u_j)\} \subset (\langle \alpha^a(u_1),\dots,\alpha^a(u_{j-1})\rangle : \alpha^a(u_j)).$
For any $k \in A(\alpha^a,u_j)$, there is $0 \leq l \leq d-1$ such that
$$i_l +a_l \leq k < i_{l+1} + a_l.$$
Let $k= k' + a_l$.
Then 
$$\frac {x_k} {x_{i_{l+1}+a_l } } \alpha^a(u_j) 
= x_{i_1} x_{i_2+a_1} \cdots x_{i_l+a_{l-1} }x_{k'+a_l}x_{i_{l+2}+a_{l+1}}\cdots x_{i_d+a_{d-1}}.$$
Since $i_l \leq k' < i_{l+1}$, it follows that $v=x_{i_1}x_{i_2}\cdots x_{i_l}x_{k'}x_{i_{l+2}}\cdots x_{i_d}$ satisfies
$v \prec u_j$  and $\alpha^a (v) =\frac {x_k} {x_{i_{l+1}+a_l } } \alpha^a(u_j)$.

On the other hand, since $I$ is strongly stable,
for any monomial $w=x_{i_1'}x_{i_2'}\cdots x_{i_d'} \in \setmono$ with $w \prec u_j$ and $w \ne u_j$,
there is $u_t\in G(I)$ such that $u_t$ divides $w$.
In particular, since $I$ is strongly stable,
we may assume $u_t=x_{i_1'}x_{i_2'}\cdots x_{i_s'} $ for some $s \leq d$.
Since $w \prec u_j$, we have $i_l' \leq i_l$ for all $1 \leq l \leq d$.
Also, since $u_t \in G(I)$ and $u_t \ne u_j$,
it follows that $u_t$ does not divide  $u_j$.
Thus $u_t$ satisfies $u_t > _{\mathrm{lex}} u_j$ and $\alpha^a(u_t)$ divides $\alpha^a(w)$.
In particular,
we have $\alpha^a(w) \in \langle \alpha^a(u_1),\dots,\alpha^a(u_{j-1}) \rangle$.

Since $v \prec u_j$ and $v \ne u_j$,
the above fact implies 
$$\alpha^a(v) \in \langle \alpha^a(u_1),\alpha^a(u_2),\dots,\alpha^a(u_{j-1}) \rangle.$$
Since $x_k \alpha^a(u_j)=x_{i_{l+1}+a_l} \alpha^a(v)$,
we have $x_k \in (\langle \alpha^a(u_1),\dots,\alpha^a(u_{j-1}) \rangle : \alpha^a(u_j)).$
\bigskip

Second, we will show $\{x_k : k\in A(\alpha^a,u_j)\}$ is a generating set of the colon ideal
$(\langle\alpha^a(u_1),\dots,\alpha^a(u_{j-1})\rangle: \alpha^a(u_j)).$
Let $w$ be a monomial belonging to the ideal $(\langle\alpha^a(u_1),\dots,\alpha^a(u_{j-1})\rangle: \alpha^a(u_j)).$
Then there is  an integer $1 \leq p <j$ such that
$\alpha^a(u_p) $ divides $w \alpha^a(u_j)$.
What we must prove is that there is $k\in A(\alpha^a,u_j)$ such that $x_k$ divides $w$.

Let $u_p =x_{j_1}x_{j_2} \cdots x_{j_{d'}}$
with $j_1 \leq j_2 \leq \cdots \leq j_{d'}$.
Since $u_p >_{\mathrm{lex}} u_j$,
there is $1 \leq r \leq d'$ such that $i_t=j_t$ for $1\leq t <r$, and $j_r <i_r$.
Since $a$ is a nondecreasing sequence,
we have $j_r+a_{r-1}<i_r + a_{r-1} \leq i_{r+1}+a_r\leq  \cdots \leq i_d + a_{d-1}$.
Then, since $\alpha^a(u_p)$ divides $w \alpha^a(u_j)$, $x_{j_r +a_{r-1}}$ must divide $w$.
On the other hand,
since $i_{r-1}=j_{r-1}\leq j_r < i_r$,
we have $$i_{r-1}+a_{r-1} \leq j_r + a_{r-1} < i_r +a_{r-1}.$$
Then $j_r+a_{r-1} \in A(\alpha^a,u_j)$ and $x_{j_r -a_{r-1}}$ divides $w$.
Hence $\{ x_k: k \in A(\alpha^a,u_j)\}$ is a generating set of
$(\langle\alpha^a(u_1),\dots,\alpha^a(u_{j-1})\rangle:\alpha^a (u_j))$.
\end{proof}

\begin{proposition} \label{satisfy}
With the same notation as in Lemma \ref{include}.
The map $\alpha^a:\setmono \to \setmono$ is a stable operator.
\end{proposition}

\begin{proof}
Let $I \subset \aru \infty$ be a finitely generated strongly stable ideal.
Then Lemma \ref{quotients} says that $\alpha^a(I)$ has linear quotients with
$$|\mathrm{set}( \alpha^a(u))|= \sum_{l=0}^{d-1}\{ i_{l+1}-i_l\}=i_d -i_0={m}(u)-1$$
for all $u=x_{i_1}x_{i_2} \cdots x_{i_d} \in G(I)$ with $i_1 \leq i_2 \leq \cdots \leq i_d$.
Then Lemma \ref{e-k} together with (\ref{hosilq}) implies
$\beta_{ij}(I) =\beta_{ij}(\alpha^a(I))$ for all $i,j$.

Also, if $J\subset I$ are finitely generated strongly stable ideals of $\aru \infty$,
then Lemma \ref{include} says
$\alpha^a (J) \subset \alpha^a (I)$.
\end{proof}

\begin{exam}
Consider the strongly stable ideal
$$I=\langle x_1^4,x_1^3x_2,x_1^3x_3,x_1^2x_2^2,x_1^2x_2x_3,x_1x_2^3,x_2^4 \rangle.$$
Let $a_1=(0,2,4,6,8,\cdots)$ and $a_2=(0,1,2,2,2,\cdots)$.
Then
$$\alpha^{a_1}(I)=
\langle x_1x_3x_5x_7,x_1x_3x_5x_8,x_1x_3x_5x_{9},x_1x_3x_6x_8,x_1x_3x_6x_9,
x_1x_4x_6x_8,x_2x_4x_6x_8 \rangle$$
and
$$\alpha^{a_2}(I)=\langle
x_1x_2x_3^2,x_1x_2x_3x_4,x_1x_2x_3x_5,x_1x_2x_4^2,x_1x_2x_4x_5,
x_1x_3x_4^2,x_2x_3x_4^2\rangle.$$

Proposition \ref{satisfy} together with Lemma \ref{e-k} says that the minimal graded free resolution of each ideal is of the form
\[
0 \longrightarrow 
 \aru n (-6)^2
\longrightarrow
 \aru n (-5)^{8}
\longrightarrow
 \aru n (-4)^{7}
\longrightarrow I
\longrightarrow 0.
\]
\end{exam}

\begin{exam}
Let $S=K[x_{ij}]_{i,j\geq 1}$ be the polynomial ring in variables $x_{ij}$ with $i \geq 1$ and $j \geq 1$.
Then $S$ is isomorphic to $\aru \infty$.
Let $p:\setmono \to S$ be the map defined by
$$p(x_1^{a_1}x_2^{a_2}\cdots x_n^{a_n})=
\prod_{i=1}^n \bigg(x_{i1}x_{i2}\cdots x_{ia_i}\bigg)$$
for any monomial $x_1^{a_1}x_2^{a_2}\cdots x_n^{a_n} \in \setmono$.
For any monomial ideal $I$,
the ideal $p(I)$ is called the \textit{polarization} of $I$.
Also, for any finitely generated monomial ideal $I$ of $\aru \infty$,
it is known that $ \beta_{ij} (p(I)) = \beta_{ij}(I)$ for all $i,j$.
Thus the polarization map $p:\setmono \to \setmono$ is a stable operator.
\end{exam}

\begin{exam}
We recall the map defined in \cite[\S 3]{Gasha}.
Let $a=(a_1,a_2,\dots)$ be a nondecreasing infinite sequence of positive integers.
Set $q_i = a_i-1$ for each $i \geq 1$.
Let $m=x_1^{\alpha_1}x_2^{\alpha_2}\cdots x_n^{\alpha_n}\in \setmono $ be a monomial of $\aru \infty$.
Then there exist integers $p_0=0<p_1<\cdots <p_c$ such that
\begin{eqnarray*}
&&\alpha_1 = q_1 + \cdots + q_{p_1 -1} +s_1  \quad\quad  \hspace{17pt}   \mbox{ with } 0\leq s_1 < q_{p_1}, \\
&&\alpha_1 = q_{p_1 +1} + \cdots + q_{p_2 -1} +s_2    \hspace{26pt}  \mbox{ with } 0\leq s_2 < q_{p_2}, \\
&& \quad \vdots \\
&&\alpha_1 = q_{p_{c-1}+1} + \cdots + q_{p_c -1} +s_c   \hspace{17pt}    \mbox{ with } 0\leq s_c < q_{p_c}.
\end{eqnarray*} 
Let
\[
 \sigma^a (x_1^{\alpha_1}x_2^{\alpha_2}\cdots x_n^{\alpha_n}) = \prod_{i=1}^c \bigg[\bigg( \prod _{j=p_{i-1}+1}^{p_i-1}x_j^{a_j}\bigg)x_{p_i}^{s_i}\bigg].
 \]
Then any monomial $ x_1^{b_1}x_2^{b_2}\cdots x_n^{b_n} \in \sigma^a(\setmono)$ satisfies $b_i <a_i$ for all $i \geq 1$.
For example,
if $I=\langle x_1^4,x_1^3x_2,x_1^2x_2^2 \rangle$ and $a=(3,3,3,\dots)$, then
$$\sigma^a(I)=\langle x_1^2x_2^2,x_1^2x_2x_3,x_1^2x_3^2 \rangle .$$
If $I$ is a finitely generated strongly stable ideal of $\aru \infty$,
then $\beta_{ij} (\sigma^a(I)) =\beta_{ij}(I)$ for all $i,j$ (\cite[Theorem 3.3]{Gasha}).
This map $\sigma^a: \setmono \to \setmono$ is also a stable operator.
\end{exam}

\section{squeezed balls and squeezed spheres}

Let $\Gamma$ be a simplicial complex on $[n]=\{1,\dots,n\}$.
Thus ${\Gamma}$
is a collection of subsets of $[n]$ such that
(i) $\{ j \} \in {\Gamma}$ for all $j \in [n]$ and
(ii) if  $T \in {\Gamma}$ and  $S \subset T$ then $S \in {\Gamma}$.
A {\em face} of $\Gamma$ is an element $S \in {\Gamma}$.
The maximal faces of $\Gamma$ under inclusion are called \textit{facets} of $\Gamma$.
A simplicial complex $\Gamma$ is called \textit{pure} if each facet of $\Gamma$ has the same cardinality.
The dimension of $\Gamma$ is the maximal integer $|S|-1$ with $S \in \Gamma$.
The \textit{$f$-vector} of a $(d-1)$-dimensional simplicial complex $\Gamma$ is the vector
$f(\Gamma)=(f_0(\Gamma),f_1 (\Gamma),\dots,f_{d-1}(\Gamma))$,
where each $f_{k-1}(\Gamma)$ is the number of faces $S$ of $\Gamma$ with $|S|=k$.

Kalai introduced squeezed spheres by extending the construction of Billera--Lee polytopes \cite{BL}.
In this section, we recall squeezed spheres.
Instead of Kalai's original definition, we use the idea which appears in \cite[\S 5.2]{K3} and \cite[pp. 246--247]{BL}.

Fix integers $d >0$ and $m\geq 0$.
Set $n=m+d+1$.
A set $U \subset \aru m$ of monomials in $\aru m$ is called an
{\em order ideal of monomials} if $U$ satisfies
\begin{itemize}
\item[(i)]
$\{1,x_1,x_2,\dots,x_m\}\subset  U$;
\item[(ii)]
if $u \in U$ and $v \in \aru m$ divides $u$, then $v \in U$.
\end{itemize}
where we let $\aru 0 =K$.
An order ideal $U \subset \aru m$ of monomials is called \textit{shifted}
if $u\in U$ and $u \prec v$ imply $v \in U$.

Let $ U \subset \aru m$ be a shifted order ideal of monomials of degree at most $\lfloor \frac {d+1} 2 \rfloor$,
where $\lfloor \frac {d+1} 2 \rfloor $ means the integer part of $\frac {d+1} 2$.
For each $u=x_{i_1}x_{i_2}\cdots x_{i_k} \in U$ with $i_1 \leq i_2 \leq \cdots \leq i_k$,
define a $(d+1)$-subset  $F_d(u) \subset [n]$ by
\begin{eqnarray*}
F_d(u)&=&
\{i_1,i_1+1\} \cup \{i_2 +2, i_2 +3\} \cup \cdots \cup \{i_k+2(k-1),i_k+2k-1\}\\
&&\cup \{n+2k-d,n+2k-d+1,\dots,n\},
\end{eqnarray*}
where $F_d(1)= \{n-d,n-d+1,\dots,n\}$.
Let $B_d(U)$ be the simplicial complex generated by
$F_d(U)=\{ F_d(u): u\in U\}.$
Kalai proved that $B_d(U)$ is a shellable $d$-ball on $[n]$.
Thus its boundary $S_d(U)=\partial (B_d(U))$ is a simplicial $(d-1)$-sphere.
This $B_d(U)$ is called a {\em squeezed $d$-ball} and $S_d(U)$ is called a {\em squeezed $(d-1)$-sphere}.
Note that every squeezed sphere is shellable (Lee \cite{L}).

We remark the reason why we assume $\{1,x_1,\dots,x_m\} \subset U$.
If $U$ is a shifted order ideal of monomials with $\{1,x_1,\dots,x_m\} \subset U$ and $x_{m+1} \not \in U$,
then $B_d(U)$ and $S_d(U)$ is the simplicial complex on the vertex set $[m+d+1]$ for $d>1$.
(That is, $\{i \} \in B_d(U) \cap S_d(U)$ for all $i=1,2,\dots,m+d+1$.)
This fact says that $m$ and $d$ determine the numbers of vertices of a squeezed sphere $S_d(U)$.
Thus, to fix the vertex set of $S_d(U)$,
we require the assumption $\{1,x_1,\dots,x_m\} \subset U$.
In particular, we will often assume $m=n-d-1$ and $S_d(U)$ is a simplicial complex on $[n]$.

We can easily know the $f$-vector of each squeezed sphere $S_d(U)$
by using $U \subset \aru m$. 
In particular $f$-vectors of squeezed spheres satisfy
the conditions of McMullen's $g$-conjecture.
To discuss $f$-vectors of squeezed spheres,
we recall $h$-vectors and $g$-vectors.

The {\em $h$-vector} $h(\Gamma)=(h_0(\Gamma),h_1(\Gamma),\dots,h_d(\Gamma))$ of a $(d-1)$-dimensional
simplicial complex $\Gamma$ is  defined by the relation
$$h_i (\Gamma) =\sum_{j=0}^i (-1)^{i-j} { d-j \choose d-i} f_{j-1}(\Gamma)
\mbox{ and } f_{i-1}(\Gamma) = \sum_{j=0}^i { d-j \choose d-i} h_i (\Gamma),$$
where we set $f_{-1}(\Gamma) =1$.
Thus, in particular, knowing the $f$-vector of $\Gamma$ is equivalent to knowing the $h$-vector of $\Gamma$.

If $\Gamma$ is a $(d-1)$-dimensional Gorenstein* complex,
then the $h$-vector of $\Gamma$ satisfies $h_i(\Gamma)=h_{d-i}(\Gamma)$ for all $0 \leq i \leq d$.
(Dehn--Sommerville equations.)
Define the {\em $g$-vector} $g(\Gamma)=(g_0(\Gamma),g_1(\Gamma),\dots,g_{ \lfloor \frac{d} {2} \rfloor}(\Gamma))$ of $\Gamma$
by
$$ g_i(\Gamma) = h_i(\Gamma)-h_{i-1}(\Gamma) \hspace{15pt}\mbox{ for } 1\leq i \leq \lfloor \frac d 2 \rfloor$$
and $g_0(\Gamma)=1$. 
Then Dehn--Sommerville equations say that
if $\Gamma$ is a Gorenstein* complex,
then knowing the $g$-vector of $\Gamma$ is equivalent to knowing the $h$-vector and the $f$-vector of $\Gamma$.

\begin{lemma}[{\cite[Proposition 5.2]{K3}}] \label{fosb}
Let $d >0$ be a positive integer,
$B_d(U)$ a squeezed $d$-ball and $S_d(U)$ a squeezed $(d-1)$-sphere.
Then
\begin{itemize}
\item[(i)]
$h_i(B_d(U))=|\{u \in U: \deg(u)=i\}|$  for all $0 \leq i \leq d+1$;
\item[(ii)]
$g_i(S_d(U))=h_i(B_d(U))=|\{u \in U : \deg(u)=i\}|$ for all $0\leq i \leq \lfloor \frac d 2 \rfloor$.
\end{itemize}
\end{lemma}

Also, the next lemma easily follows.

\begin{lemma} \label{skeleton}
Let $m\geq 0$, $d > 0$ and $ \lfloor \frac{d+1} 2 \rfloor \geq k \geq 0$ be integers.
Let $U \subset \aru m$ be a shifted order ideal of monomials of degree at most $k$.
Then we have $f_{i-1}(B_d(U))=f_{i-1}(S_d(U))$ for $1 \leq i \leq d-k$.
\end{lemma}

\begin{proof}
Lemma \ref{fosb} together with Dehn--Sommerville equations of $S_d(U)$ says
\begin{eqnarray*}
h_i(S_d(U)) -h_{i-1}(S_d(U))=
\left\{
\begin{array}{ll}
h_i( B_d(U)), \hspace{70pt} \mbox{ for } 0\leq i \leq \lfloor \frac d 2 \rfloor,\\
0, \hspace{116pt} \mbox{ for } i= \frac {d+1} 2,\\
-h_{d+1-i}(B_d(U)), \hspace{39pt} \mbox{ for } \lfloor \frac {d+1} 2 \rfloor  <i \leq d,
\end{array}
\right.
\end{eqnarray*}
where we let $h_{-1}(S_d(U))=0$.
Lemma \ref{fosb} also says 
$$h_i(B_d(U)) = - h_{d+1-i}(B_d(U))=0 \mbox{\ \ \  for } \ k +1 \leq i \leq d-k.$$
Thus we have 
$$h_i(S_d(U)) -h_{i-1}(S_d(U))= h_i(B_d(U)) \mbox{ \ \ \ for } 0 \leq i \leq d-k.$$
Since $ { {d-j} \choose {d-i}} = {{d-j+1} \choose {d-i+1}} -{ {d-j} \choose {d-i+1}}$,
we have 
\begin{eqnarray*}
f_{i-1} (S_d(U))&=& \sum_{j=0}^i { {d-j} \choose {d-i}} h_j(S_d(U)) \\
&=&  \sum_{j=0}^i \bigg\{ { {d-j+1} \choose {d-i+1}} - {{d-j} \choose {d-i+1}} \bigg\}  h_j(S_d(U)) \\
&=& \sum_{j=0}^i { {d-j+1} \choose {d-i+1} }  \{h_j(S_d(U))-h_{j-1} (S_d(U))\}\\
&=& f_{i-1}(B_d(U))
\end{eqnarray*}
for $1 \leq i \leq d-k$, as desired.
\end{proof}

\section{Lefschetz properties}

In this section,
we recall some facts about generic initial ideals and Lefschetz properties.
We refer the reader to \cite{St} for the fundamental theory of Stanley--Reisner rings, Cohen--Macaulay complexes and
Gorenstein* complexes.

The \textit{Stanley--Reisner ideal} $I_\Gamma \subset \aru n$ 
of a simplicial complex $\Gamma$ on $[n]$ is a monomial ideal generated by
all squarefree monomials $x_{i_1}x_{i_2}\cdots x_{i_k} \in \aru n$ with $\{i_1,i_2,\dots,i_k\} \not \in \Gamma$.
The quotient ring $R(\Gamma)= \aru n / I_\Gamma$ is called the \textit{Stanley--Reisner ring of $\Gamma$}.



Let $\Gamma$ be a $(d-1)$-dimensional Gorenstein* complex on $[n]$.
We say that $\Gamma$ has the {\em weak Lefschetz property}
if there is a system $\vartheta_1,\dots,\vartheta_d$ of parameters
of $\aru n / I_\Gamma$ and a linear form $\omega \in \aru n$ such that the multiplication map 
$\omega: H_{i-1}(\Gamma) \to H_i(\Gamma)$ is injective for $1\leq i \leq \lfloor \frac d 2 \rfloor$ 
and surjective for $i > \lfloor \frac d 2 \rfloor$,
where $H_i(\Gamma)$
is the $i$-th homogeneous component of $(\aru n /I_\Gamma)\otimes (\aru n /\langle\vartheta_1,\dots,\vartheta_d\rangle)$.
Also, we say that $\Gamma$ has the {\em strong Lefschetz property}
if there is a system $\vartheta_1,\dots,\vartheta_d$ of parameters of $\aru n / I_\Gamma$ and
there is a linear form $\omega \in \aru n$ such that the multiplication map $\omega^{d-2i}:H_i(\Gamma) \to H_{d-i} (\Gamma)$ is an isomorphism for $ 0 \leq i \leq  
\lfloor \frac d 2\rfloor$.
The above linear form $\omega$ is called a\textit{ weak (resp. strong) Lefschetz element} of
$(\aru n / I_\Gamma) \otimes (\aru n /\langle\vartheta_1,\dots,\vartheta_d\rangle)$.
The following fact is well known.

\begin{lemma} \label{glf}
Let $\Gamma$ be a $(d-1)$-dimensional Gorenstein* complex  on $[n]$ with the weak (resp. strong) Lefschetz property
and $\vartheta_1,\vartheta_2,\dots,\vartheta_d,\omega$ generic linear forms of $\aru n$.
Then $\vartheta_1,\vartheta_2,\dots,\vartheta_d$ is a system of parameters of $R(\Gamma)$
and $\omega$ is a weak (resp. strong) Lefschetz element of $(\aru n / I_\Gamma) \otimes (\aru n /\langle\vartheta_1,\dots,\vartheta_d \rangle)$.
\end{lemma}


For a $(d-1)$-dimensional simplicial complex $\Gamma$ on $[n]$,
define the set of monomials $L(\Gamma) = \bigcup _{i \geq 0}L_i(\Gamma)$ by
$$ L_i(\Gamma)= \{u \in \aru {n-d}:u\not\in \gin(I_\Gamma) \mbox{ is a monomial of degree }i\}$$
and
define  $U(\Gamma) = \bigcup_{i\geq 0}U_i (\Gamma)$ by
$$U_i(\Gamma)= \{ u \in \aru {n-d-1}:u \not \in \gin(I_\Gamma) \mbox{ is a monomial of degree } i\}.$$

\begin{lemma} \label{hcm}
Let $\Gamma$ be a $(d-1)$-dimensional Cohen--Macaulay complex on $[n]$.
Then
\begin{itemize} 
\item[(i)]$\max(\gin (I_\Gamma)) = {n-d}$;
\item[(ii)]
$|L_i(\Gamma)|=h_i(\Gamma)$ for $i\geq 0$, where we let $h_i(\Gamma) =0$ for $i > d$.
\end{itemize}
\end{lemma}

\begin{proof}
(i) Since $\Gamma$ is a Cohen--Macaulay complex,
the Auslander--Buchsbaum formula says $\mathrm{proj \ dim}(\aru n / I_\Gamma) = n-d$.
Then Lemma \ref{idealmax} says $\max(\gin(I_\Gamma)) = n-d$.
\medskip

(ii)
Let $\vartheta_1,\vartheta_2,\dots,\vartheta_d$ be a system of parameters of $\aru n / I_\Gamma$ which are linear forms.
Since $\max(\gin(I_\Gamma)) \leq n-d$,
the sequence $x_{n-d+1},x_{n-d+2},\dots,x_n$ is a system of parameters of $\aru n / \gin(I_\Gamma)$.
Since $\aru n / I_\Gamma$ and $\aru n / \gin(I_\Gamma)$ have the same Hilbert function, 
$(\aru n / I_\Gamma) \otimes (\aru n /\langle \vartheta_1,\dots,\vartheta_d \rangle)$ and
$(\aru n / \gin(I_\Gamma)) \otimes (\aru n /\langle x_{n-d+1},\dots,x_n \rangle)$ have the same Hilbert function.
On the other hand, it is well known \cite[pp. 57--58]{St} that
$h_i(\Gamma)$ is equal to the $i$-th Hilbert function of  $(\aru n / I_\Gamma) \otimes (\aru n /\langle \vartheta_1,\dots,\vartheta_d \rangle)$.
Since $L_i(\Gamma)$ is a $K$-basis of the $i$-th homogeneous component of $(\aru n / \gin(I_\Gamma)) \otimes (\aru n/\langle x_{n-d+1},\dots,x_n \rangle)$,
we have $|L_i(\Gamma) | = h_i(\Gamma)$ for all $i \geq 0$.
\end{proof}

Lemma \ref{hcm} (i) says that if $\Gamma$ is Cohen--Macaulay,
then $\gin(I_\Gamma) \cap \aru {n-d}$ (or $L(\Gamma)$)
determines $\gin(I_\Gamma)$.
The next lemma immediately follows from \cite[Lemma 2.7 and Proposition 2.8]{W}.

\begin{lemma} \label{lefschetz1}
Let $\Gamma$ be a $(d-1)$-dimensional Gorenstein* complex on $[n]$.
For any integer $k\geq 1$, we write $x_{n-d}^kL_{i}(\Gamma)=\{ x_{n-d}^k u: u \in L_i(\Gamma)\}$.
Then
\begin{itemize}
\item[(i)]
$\Gamma$ has the weak Lefschetz property if and only if 
$x_{n-d}L_{i-1}(\Gamma) \subset L_i(\Gamma)$ for $1\leq i \leq \lfloor \frac d 2 \rfloor$
and
$x_{n-d}L_{i-1}(\Gamma) \supset L_i(\Gamma)$  for $i > \lfloor \frac d 2 \rfloor$.
\item[(ii)]
$\Gamma$ has the strong Lefschetz property if and only if 
$x_{n-d}^{d-2i}L_i(\Gamma) = L_{d-i}(\Gamma)$ for $0 \leq i \leq \lfloor \frac d 2 \rfloor$.
\end{itemize}
\end{lemma}

\begin{lemma} \label{lefschetz2}
Let $\Gamma$ be a $(d-1)$-dimensional Gorenstein* complex on $[n]$.
Then
\begin{itemize}
\item[(i)]
$\Gamma$ has the weak Lefschetz property if and only if
\begin{eqnarray*}
|U_i(\Gamma)|=
\left\{
\begin{array}{l}
g_i(\Gamma),
\hspace{20pt}\mbox{for } 0 \leq i \leq \lfloor \frac d 2 \rfloor, \\
\hspace{0cm} {0}, \hspace{40pt}
\mbox{for } i > \lfloor \frac d 2 \rfloor.
\end{array}
\right.
\end{eqnarray*}
\item[(ii)]
Assume that $\Gamma$ has the strong Lefschetz property.
Let $u=u'x_{n-d}^t \in \aru {n-d}$ be a monomial with $\deg(u)=k$.
Then $u \in L(\Gamma)$ if and only if $u' \in U(\Gamma)$ and $t \geq 2k-d$.
\end{itemize}
\end{lemma}

\begin{proof}(i)
Since $L(\Gamma)$ is an order ideal of monomials,
if a monomial $x_{n-d}u \in L(\Gamma)$ then $u\in L(\Gamma)$.
This fact says $L_i(\Gamma) \subset U_i(\Gamma) \cup x_{n-d}L_{i-1}(\Gamma)$.
Since $U_i(\Gamma) \cap x_{n-d}L_{i-1}(\Gamma)= \emptyset$,
Lemma \ref{hcm} says
$$|U_i(\Gamma)|\geq|L_i(\Gamma)|-|L_{i-1}(\Gamma)|=h_{i}(\Gamma)-h_{i-1}(\Gamma).$$
Then it follows that $|U_i(\Gamma)|=h_i(\Gamma)-h_{i-1}(\Gamma)$ if and only if $L_i(\Gamma) \supset x_{n-d}L_{i-1}(\Gamma)$.

Also, it is easy to see that
$|U_i(\Gamma)|=0$ if and only if $L_i(\Gamma) \subset x_{n-d}L_{i-1}(\Gamma)$.
Thus the assertion follows from Lemma \ref{lefschetz1}.
\medskip

(ii) 
Assume $k \leq \lfloor \frac d 2 \rfloor$.
Since $\Gamma$ has the weak Lefschetz property,
Lemma \ref{lefschetz1} together with the proof of (i) says
$$L_k(\Gamma)= U_k(\Gamma) \cup x_{n-d}L_{k-1}(\Gamma).$$ 
Thus, $u\in L_k(\Gamma)$ if and only if $u\in U_k(\Gamma)$ or $u/x_{n-d} \in L_{k-1}(\Gamma)$.
Inductively, we have $u \in L_k(\Gamma)$ if and only if $u' \in U(\Gamma)$.

Assume $k > \lfloor \frac d 2 \rfloor$.
Since $\Gamma$ has the strong Lefschetz property, we have $L_k(\Gamma) =x_{n-d}^{2k-d}L_{d-k} (\Gamma)$.
Thus $u \in L_k(\Gamma)$ if and only if $t \geq 2k-d$ and $u/x_{n-d}^{2k-d} \in L_{d-k}(\Gamma)$.
Since $d-k \leq \lfloor \frac d 2 \rfloor$,
we have $u/x_{n-d}^{2k-d} \in L_{d-k}(\Gamma)$ if and only if $u' \in U(\Gamma)$.
\end{proof}

Lemma \ref{lefschetz2} says that if $\Gamma$ is a Gorenstein* complex with the strong Lefschetz property
then $U(\Gamma)$ determines $\gin(I_\Gamma)$.

\section{generic initial ideals of squeezed balls and squeezed spheres}

If $\Gamma$ is a $(d-1)$-dimensional Gorenstein* complex on $[n]$ with the weak Lefschetz property
and if $K$ is a field of characteristic $0$,
then Lemma \ref{lefschetz2} says
$$U(\Gamma)= \{  u\in \aru {n-d-1}: u \not\in \gin(I_\Gamma)  \mbox{ is a monomial}\}\subset \aru {n-d-1}$$
is a shifted order ideal of monomials of degree at most $\lfloor \frac d 2 \rfloor$
with $g_i(\Gamma)=|\{u\in U(\Gamma): \deg(u)=i\}|$ for  $0 \leq i \leq \lfloor \frac d 2 \rfloor$.
(Since $g_1 (\Gamma)=n-d-1$ and $g_0 (\Gamma) =1$, the set $U(\Gamma)$ certainly contains $\{1,x_1,\dots,x_{n-d-1}\}$.)
Furthermore, if $\Gamma$ has the strong Lefschetz property,
then Lemma \ref{lefschetz2} also says that $U(\Gamma)$ determines $\gin(I_\Gamma)$.

Conversely,
for any shifted order ideal $U \subset \aru {n-d-1}$ of monomials of degree at most $\lfloor \frac d 2 \rfloor$,
there is the squeezed $(d-1)$-sphere $S_d(U)$ on $[n]$ with 
$g_i(S_d(U))=|\{u\in U: \deg (u)=i\}|$ for $0 \leq i \leq \lfloor \frac d 2 \rfloor$.
Kalai conjectured that
\[ U(S_d(U))=U. \]

In this section, we will prove this equality.
Fix positive integers $n>d>0$.
Let $ U\subset \aru {n-d-1}$ be a shifted order ideal of monomials of degree at most $\lfloor \frac {d+1} 2 \rfloor$.
Write $I(U)\subset \aru n$ for the ideal generated by all monomials $u \in \aru {n-d-1}$ with $u \not \in U$.
Since $U$ is shifted, $I(U)$ is a strongly stable ideal.

Let $\alpha^2:\setmono \to \setmono$ be the map defined by
$$\alpha^2(x_{i_1}x_{i_2}x_{i_3} \cdots x_{i_k})=x_{i_1}x_{i_2+2}x_{i_3+4}\cdots x_{i_k + 2(k-1)},$$
for any $x_{i_1}x_{i_2}\cdots x_{i_k} \in \setmono$ with $i_1 \leq i_2 \leq \dots \leq i_k$.
Then, by Proposition \ref{satisfy}, the map $\alpha^2:\setmono \to \setmono$ is a stable operator.
Since $\max\{ \deg(u) :u \in G(I(U))\} \leq \lfloor \frac{d+1} 2\rfloor +1$,
we have $m(\alpha^2(u)) \leq n-d-1+ 2( \lfloor \frac {d+1} 2 \rfloor )\leq n$ for all $u \in G(I(U))$.
We write $\alpha^2(I(U))$ for the ideal of $\aru n$ generated by $\{ \alpha^2(u) : u\in G(I(U))\}$.

\begin{proposition} \label{ginball}
Let $n>d>0 $ be positive integers.
Let $B_d(U)$ be a squeezed $d$-ball on $[n]$ and $I(U) \subset \aru n$ 
 the ideal generated by all monomials $u \in \aru {n-d-1}$ with $u \not \in U$.
Then one has
$$I_{B_d(U)}=\alpha^2(I(U)).$$
In particular, one has $\gin(I_{B_d(U)})=I(U)$.
\end{proposition}

\begin{proof}
First, we will show $I_{B_d(U)}\supset \alpha^2(I(U))$.
Let $u=x_{i_1}x_{i_2}\cdots x_{i_k}\in G(I(U))$ 
with $i_1 \leq i_2 \leq \dots \leq i_k$
and let 
\begin{eqnarray}
f(u)=\{i_1,i_2+2,i_3+4,\dots,i_k+2(k-1)\}. \label{kyuu}
\end{eqnarray}
Note that
\begin{eqnarray}
\max (f(u))\leq n-d-1+2(k-1) =n+2k-d-3. \label{po1}
\end{eqnarray}
We will show $f(u) \not\in B_d(U)$.
Since $B_d(U)$ is the simplicial complex generated by $\{ F_d(v):v\in U\}$,
if $f(u) \in B_d(U)$ then there is $ w \in U$ such that
$f(u) \subset F_d(w)$.
Thus what we must prove is  $f(u) \not \subset F_d(w)$ for all $w \in U$.

Let $w=x_{j_1}x_{j_2}\cdots x_{j_l}\in U$ with $j_1 \leq j_2 \leq \cdots \leq j_l$.
Set 
$$F^\flat (w)=\{j_1,j_1+1\}\cup \{j_2+2,j_2+3\} \cup \{j_3+4,j_3+5\}\cup\cdots\cup \{j_l+2(l-1),j_l+2l-1\}$$
and
$$F^\natural(w)=\{n+2l-d,n+2l-d+1,\dots,n\}.$$
Then $F_d(w)=F^\flat(w) \cup F^\natural(w).$
\medskip

[Case 1]
Assume  $l <k$.
Then the form (\ref{kyuu}) says $|F^\flat(w)\cap f(u)|\leq l$.
On the other hand,
(\ref{po1}) says
\begin{eqnarray*}
|F^\natural(w) \cap f(u)| &=& |\{n+2l-d,n+2l-d+1,\dots,n+2k-d-3\}\cap f(u)|.
\end{eqnarray*}
Since $|\{ n+2l-d,n+2l-d+1,\dots,n+2k-d-3\}| =2(k-l-1)$,
the form (\ref{kyuu}) says 
$|F^\natural(w) \cap f(u)| \leq k-l-1.$
Thus we have $|F_d(w) \cap f(u)| \leq k-1$ and $f(u) \not\subset F_d(w)$.
\medskip

[Case 2]
Assume $l \geq k$ and $f(u) \subset F_d(w)$.
Then, the form (\ref{kyuu}) says
$$j_p +2(p-1) \leq i_p +2(p-1)$$
for all $1 \leq p \leq k.$
Thus we have $x_{j_1}x_{j_2} \cdots x_{j_k} \preceq u$. 
Since $U$ is shifted and since $u \not\in U$, we have $x_{j_1}x_{j_2}\cdots x_{j_k} \not \in U$.
However, since $U$ is an order ideal of monomials and since $w \in U$,
$x_{j_1}x_{j_2} \cdots x_{j_k}$ must be contained in $U$.
This is a contradiction.
Thus we have $f(u) \not \subset F_d(w)$.
\medskip

Second,
we will show that $I_{B_d(U)}$ and $\alpha^2 (I(U))$ have the same Hilbert function.

Lemma \ref{hcm} says that the $i$-th Hilbert function of $\aru {n-d-1} / (\gin(I_{B_d(U)}) \cap \aru {n-d-1})$
is equal to $h_i(B_d(U))$ for all $i \geq 0$,
where we let $h_i(B_d(U))=0$ for $i > d+1$.
On the other hand, the $i$-th Hilbert function of $\aru {n-d-1} / (I(U)) \cap \aru {n-d-1})$
is equal to $|\{u\in U: \deg(u)=i\}|$.
Then Lemma \ref{fosb} says $I(U) \cap \aru {n-d-1}$ and $\gin (I_{B_d(U)}) \cap \aru {n-d-1}$ have 
the same Hilbert function.
Since $\max(I(U)) \leq {n-d-1}$ by the definition of $I(U)$ and $\max( \gin (I_{B_d(U)})) \leq n-d-1$
by Lemma \ref{hcm},
it follows that $I(U)$ and $\gin(I_{B_d(U)})$ have the same Hilbert function.
Since Theorem \ref{soperator} says that $\gin(\alpha^2(I(U)))=I(U)$,
it follows from Lemma \ref{bnt} that
$I_{B_d(U)}$ and $\alpha^2 (I(U))$ have the same Hilbert function.
\medskip

Then we proved that $I_{B_d(U)} \supset  \alpha^2(I(U))$ and $I_{B_d(U)}$ and $\alpha^2(I(U))$ have the same Hilbert function.
Thus we have $I_{B_d(U)} =  \alpha^2(I(U))$.
In particular, Theorem \ref{soperator} guarantees $\gin(I_{B_d(U)})=I(U)$.
\end{proof}

\begin{theorem} \label{conj}
Let $n>d>0$ be positive integers.
Let $U\subset \aru {n-d-1}$ be a shifted order ideal of monomials of degree at most $\lfloor \frac d 2 \rfloor$.
Then one has
$$U(S_d(U)) =U.$$
\end{theorem}

\begin{proof}
Since $S_d(U) \subset B_d(U)$ and since $d - \lfloor \frac d 2 \rfloor= \lfloor \frac {d+1} {2} \rfloor$,
Lemma \ref{skeleton} says
$$\{S \in S_d(U) : |S|\leq \lfloor \frac {d+1}  2 \rfloor\}
=\{S \in B_d(U):|S| \leq \lfloor \frac {d+1}  2 \rfloor\}.$$
Thus we have 
$(I_{S_d(U)}) _{\leq \lfloor \frac {d+1}  2 \rfloor} = (I_{B_d(U)})_{\leq \lfloor \frac {d+1}  2 \rfloor},$
where $I_ {\leq k}$ denotes the ideal generated by all polynomials $f$ in a graded ideal $I \subset \aru n$
with $\deg(f) \leq k$. 
Then, Proposition \ref{ginball} says
$$\gin(I_{S_d(U)})_{\leq \lfloor \frac {d+1}  2 \rfloor} =\gin(I_{B_d(U)})_{\leq \lfloor \frac {d+1} 2 \rfloor}
= I(U)_{\leq \lfloor \frac {d+1}  2 \rfloor}.$$
The construction of $I(U)$ says that $I(U)$ contains all monomials $u\in \aru {n-d-1}$
with $\deg (u) \geq \lfloor \frac d 2 \rfloor +1$.
Since $\gin(I_{S_d(U)} ) \supset \gin(I_{B_d(U)})=I(U)$,
we have
\begin{eqnarray*}
U(S_d(U))&=& \{ u \in \aru {n-d-1}: u\not\in \gin(I_{S_d(U)}) \mbox{ is a monomial}\} \\
&=& \{u \in \aru {n-d-1}: u\not\in \gin(I_{S_d(U)})\mbox{ is a monomial of degree } \leq \lfloor \frac d 2 \rfloor\} \\
&=& \{u \in \aru {n-d-1}: u\not\in I(U)\mbox{ is a monomial of degree }\leq \lfloor \frac d 2 \rfloor\} \\
&=&U,
\end{eqnarray*}
as desired.
\end{proof}

The squeezed spheres considered in Theorem \ref{conj} are in fact special instances of the general case.
This is because for general squeezed spheres, we assume that
$U \subset \aru {n-d-1}$ is a shifted order ideal of monomials of 
degree at most $\lfloor \frac{d+1} 2 \rfloor$.
We call $S_d(U)$ a {\em special squeezed $(d-1)$-sphere} (\textit{S-squeezed $(d-1)$-sphere}) if $U \subset \aru {n-d-1}$ is a shifted order ideal of
monomials of degree at most $\lfloor \frac d 2 \rfloor$.
If $d$ is even, then every squeezed $(d-1)$-sphere is an S-squeezed $(d-1)$-sphere.
Also, it is easy to see that
$S_d(U)$ is an S-squeezed $(d-1)$-sphere if and only if $B_d(U)$ is the cone over $B_{d-1}(U)$,
that is,
$B_d(U)$ is generated by $\{ \{n\} \cup F_{d-1}(u):u \in U\}$.
Theorem \ref{conj} and Lemma \ref{lefschetz2} (i) imply the following corollaries.

\begin{cor} \label{csswl}
Every S-squeezed sphere has the weak Lefschetz property.
\end{cor}

\begin{cor} \label{charawl}
Let $K$ be a field of characteristic $0$ and  $n>d>0$ positive integers.
A shifted order ideal $U\subset \aru {n-d-1}$ of monomials is equal to  $U=U(\Gamma)$ 
for some $(d-1)$-dimensional Gorenstein* complex (or for some simplicial $(d-1)$-sphere) $\Gamma$ on $[n]$ with the weak Lefschetz property
if and only if
$U$ is a shifted order ideal of monomials of degree at most $\lfloor \frac d 2 \rfloor$.
\end{cor}

Although we only proved that S-squeezed spheres have the weak Lefschetz property,
it seems likely that Corollary \ref{charawl} is true
when we consider the strong Lefschetz property instead of the weak Lefschetz property.
The remaining problem is

\begin{problem}
Prove that every squeezed sphere (or every S-squeezed sphere) has the strong Lefschetz property.
\end{problem}

If a Gorenstein* complex $\Gamma$ has the strong Lefschetz property,
then $U(\Gamma)$ determines $\gin(I_\Gamma)$.
Thus the above problem would yield a complete characterization of generic initial ideals
of Stanley--Reisner ideals of Gorenstein* complexes with the strong Lefschetz property,
when the base field is of characteristic $0$.

\section{The squeezing operation and graded Betti numbers}

Let $K$ be a field of characteristic $0$.
Let $\Gamma$ be a $(d-1)$-dimensional Gorenstein* complex on $[n]$ with the weak Lefschetz property.
Then
$U(\Gamma)= \{ u \in \aru {n-d-1}:u\not\in \gin(I_\Gamma) \mbox{ is a monomial}\}$ 
is a shifted order ideal of monomials of degree at most $\lfloor \frac d 2 \rfloor$. 
Define 
$$\mathrm{Sq}(\Gamma) = S_d(U(\Gamma)).$$
Then Lemmas \ref{fosb} and \ref{lefschetz2} say that
$\Gamma$ and $\mathrm{Sq} (\Gamma)$ have the same $f$-vector.
This operation $\Gamma \to \mathrm{Sq}(\Gamma)$ is called \textit{squeezing}.

The squeezing operation was considered by Kalai.
Since it is conjectured that every simplicial sphere has the weak Lefschetz property,
it is expected that squeezing becomes an operation for simplicial spheres
and acts like a shifting operation
(see \cite{K} for shifting operations).
In the present paper,
we study the behavior of graded Betti numbers under squeezing.



In this section,
we write $\beta_{ij}^R (M)$ for the graded Betti numbers of a graded $R$-module $M$ over a graded ring $R$.
Let $\Gamma$ be a $(d-1)$-dimensional Gorenstein* complex with the weak Lefschetz property,
$\vartheta_1,\vartheta_2,\dots,\vartheta_d$ generic linear forms
and $\bar R =  \aru n / \langle \vartheta_1,\vartheta_2,\dots,\vartheta_d \rangle$.
Let
$A= (\aru n / I_\Gamma) \otimes \bar R$.
Then $A$ is a $0$-dimensional Gorenstein ring with
$\beta_{ij}^{\aru n}(\aru n/ I_\Gamma)=\beta_{ij}^{\bar R}(A)$ for all $i,j$.
The following fact is known.

\begin{lemma}[{\cite[Proposition 8.7]{MN}}] \label{gbetti}
With the notation as above.
Let $\omega\in \aru n$ be a linear form.
Set $\tilde R= \aru n / \langle \vartheta_1,\vartheta_2,\dots,\vartheta_d,\omega \rangle$.
Then
\begin{eqnarray*}
\beta_{ii+j}^{\aru n}( \aru n / I_{\Gamma})
\leq\beta_{ii+j}^{\tilde R}(A / \omega A) +\beta_{n-d-i,n-i-j}^{\tilde R}(A / \omega A), \hspace{30pt} 
\mbox{ for all } \ i \mbox{ and }j.\\
\end{eqnarray*}
\end{lemma}

Let $U\subset \aru {n-d-1}$ be a shifted order ideal of monomials of degree at most $\lfloor \frac d 2 \rfloor$ and
$I(U)\subset \aru n$ the ideal generated by all monomials $u \in \aru {n-d-1}$ with $u\not\in U$.
Since $I(U)$ is strongly stable, we can easily compute the graded Betti numbers of $\aru n / I(U)$ by
the Eliahou--Kervaire formula \cite{EK}.

An order ideal $U \subset \aru {n-d-1}$ of monomials is called a \textit{lexicographic order ideal of monomials}
if $u \in U$ and $v <_{\mathrm{lex}} u$ imply $v \in U$ for all monomials $u$ and $v$ in $\aru {n-d-1}$ with $\deg (u) =\deg(v)$.
If $U$ is a lexicographic order ideal of monomials of degree at most $\lfloor \frac d 2 \rfloor $,
then $S_d(U)$ is the boundary complex of a simplicial $d$-polytope, called the Billera--Lee polytope \cite{BL}.
Migliore and Nagel proved that the graded Betti numbers of the Stanley--Reisner ideal of the boundary complex
$S_d(U)$ of any Billera--Lee polytope are easily computed by using $I(U)$.

\begin{lemma}[{\cite[Theorem 9.6]{MN}}] \label{blbetti}
Let $S_d(U)$ be the boundary complex of a Billera--Lee $d$-polytope on $[n]$, $R = \aru n$ and
$I(U)\subset R$ the ideal generated by  all monomials $u \in \aru {n-d-1}$ with $u \not\in U$.
Then
\begin{eqnarray*}
\beta_{ii+j}^{R}( R / I_{S_d(U)})=
\left\{
\begin{array}{l}
\beta_{ii+j}^{R}(R / I(U)), \hspace{138pt}
\mbox{for} \ j < \frac d 2, \\
\beta_{ii+j}^{R}(R/ I(U)) +\beta_{n-d-i,n-i-j}^{R}(R/ I(U)), \hspace{10pt} 
\mbox{ for } j= \frac d 2, \\
\beta_{n-d-i,n-i-j}^{R}(R / I(U)), \hspace{102pt}
\mbox{for } j > \frac d 2 .
\end{array}
\right.
\end{eqnarray*}
\end{lemma}

We will show the same property for S-squeezed spheres.
\begin{theorem} \label{sqbetti}
Let $S_d(U)$ be an S-squeezed $(d-1)$-sphere on $[n]$, $R=\aru n$ and 
$I(U)\subset R$ the ideal generated by all monomials  $u \in \aru {n-d-1}$ with $u \not\in U$.
Then
\begin{eqnarray*}
\beta_{ii+j}^{R}( R / I_{S_d(U)})=
\left\{
\begin{array}{l}
\beta_{ii+j}^{R}(R / I(U)), \hspace{138pt}
\mbox{for} \ j  <  \frac {d} 2, \\
\beta_{ii+j}^{R}(R/ I(U)) +\beta_{n-d-i,n-i-j}^{R}(R/ I(U)), \hspace{10pt} 
\mbox{ for } j= \frac d 2, \\
\beta_{n-d-i,n-i-j}^{R}(R / I(U)), \hspace{100pt}
\mbox{for } j >  \frac d 2.
\end{array}
\right.
\end{eqnarray*}
\end{theorem}

\begin{proof}
It follows from  \cite[Lemma 1.2]{HerzogHibi} that,
for any graded ideal $I \subset R$ and for any integer $k \geq 0$,
one has
$$\beta_{ii+j}^R (I)  = \beta_{ii+j}^R(I _{\leq k})\ \ \ \mbox{ for } j \leq k.$$
On the other hand, Lemma \ref{skeleton} and Proposition \ref{ginball} say
\begin{eqnarray}
(I_{S_d(U)}) _{\leq \lfloor \frac {d+1} 2 \rfloor} =
(I_{B_d(U)})_{\leq \lfloor \frac {d+1} 2 \rfloor}= \alpha^2 (I(U)) _{ \leq \lfloor \frac {d+1} 2 \rfloor}. \label{tuika}
\end{eqnarray}
Recall that, for any graded ideal $J \subset R$,
 one has $\beta_{i+1k} (R / J) =\beta_{ik}(J)$ for all $i \geq 0$.
Then, since $\alpha^2: \setmono \to \setmono$ is a stable operator,
the equality (\ref{tuika}) says
$$ \beta_{ii+j}^R(R/ I_{S_d(U)})= \beta_{ii+j}^R (R / \alpha^2(I(U))) = \beta_{ii+j}^R(R/I(U))$$
for all $j  < \frac {d} 2$.
Then, by the self duality of the Betti numbers of Gorenstein rings,
we have
\begin{eqnarray}
\beta_{ii+j}^{R}(R / I_{S_d(U)})=
\left\{
\begin{array}{l}
\beta_{ii+j}^{R}(R/ I(U)), \hspace{80pt}\mbox{ for } j < \frac {d} 2, \\
\beta_{n-d-i,n-i-j}^{R}(R/I(U)), \hspace{43pt}\mbox{ for } j > \frac {d} 2.
\end{array}
\right. \label{hosi61}
\end{eqnarray}
Thus the only remaining part is $j= \frac d 2$ when $d$ is even.
Note that this part must be determined by the Hilbert function of $R/ I_{S_d(U)}$.

We use the following well known fact:
For any graded ideal $I \subset R$,
let $a_k(R / I )= \sum_{i=0}^k (-1)^i \beta_{ik}^{R}(R/I)$ for $k \geq 0$.
Then the Hilbert function $H(R/I,t)$ of $R/I$ is given by
$H(R/I,t)= \sum _{j \geq 0} a_j(R/I) { n-1 + t -j \choose t-j}.$
Also, this $a_k(R/I)$ is uniquely determined by the Hilbert function of $R/I$.
(See \cite[Lemma 4.1.13]{CM}.)

On the other hand, for any shifted order ideal $U\subset \aru {n-d-1}$ of monomials,
there is the unique lexsegment order ideal $U^{\mathrm{lex}}\subset \aru {n-d-1}$ of monomials such that
$|\{u\in U: \deg(u)=k\}|= |\{ u\in U^{\mathrm{lex}}:\deg(u)=k\}|$ for all $k \geq 0$
(see \cite[\S 2]{BL}).
Then Lemma \ref{fosb} says that $I_{S_d(U)}$ and $I_{S_d(U^{\mathrm{lex}})}$ have the same Hilbert function.

Also, since $I(U)\cap \aru {n-d-1}$ and $I(U^{\mathrm{lex}})\cap \aru {n-d-1}$ have the same Hilbert function
and since $\max(I(U))\leq n-d-1$ and $\max(I(U^{\mathrm{lex}})) \leq n-d-1$,
it follows that $I(U)$ and $I(U^{\mathrm{lex}})$ have the same Hilbert function.
Thus we have $a_k(R/I_{S_d(U)}) = a_k(R/I_{S_d(U^{\mathrm{lex}})})$ and 
$a_k(R/I(U))= a_k (R/I(U^{\mathrm{lex}}))$ for all $k \geq 0$.

Since $I(U)$ and $I(U^{\mathrm{lex}})$ are strongly stable and
since they have no generator of degree $> \frac d 2 +1$,
the Eliahou--Kervaire formula says that $\beta_{ii+j}^R(R/I(U))=0$ and $\beta_{ii+j}^R(R/I(U^{\mathrm{lex}}))=0$ for $j \geq  \frac {d} 2 +1$.
Then Lemma \ref{blbetti} says
$$a_k(R/I(U^{\mathrm{lex}})) + (-1)^{n-d} a_{n-k}(R/I(U^{\mathrm{lex}}))=a_k(R/I_{S_d(U^{\mathrm{lex}})})$$
for all $k \geq 0$.
Thus we have \begin{eqnarray}
a_k(R/I(U)) + (-1)^{n-d} a_{n-k}(R/I(U))=a_k(R/I_{S_d(U)}) \label{aka}
\end{eqnarray}
for all $k \geq 0$.
Then, by using (\ref{hosi61}) and (\ref{aka}), a routine computation says  
$$\beta_{ii+\frac d 2}^{R}(R/I(U)) + \beta_{n-d-i, n-i-\frac d 2}^{R}(R/ I(U))
=\beta_{ii+\frac d 2}^{R} (R/ I_{S_d(U)})$$
for all $i\geq 0$, as desired.
\end{proof}

Next, we will show that squeezing increases graded Betti numbers.
Before the proof, we recall the important relation between generic initial ideals and
generic hyperplane sections.

Let $h_1=\sum_{j=1}^n a_j x_j$ be a linear form of $\aru n$ with $a_n \ne 0$.
Define a homomorphism $\Phi_{h_1} :\aru n \to \aru {n-1}$ by
$\Phi_{h_1}(x_j)=x_j$ for $1 \leq j \leq n-1$ and $\Phi_{h_1} (x_n) = - \frac 1 {a_n} ( \sum_{j=1}^{n-1} a_j x_j)$.
Then $\Phi_{h_1}$ induces a ring isomorphism between $(\aru n /\langle h_1 \rangle)$ and $\aru {n-1}$.
Let $f\in \aru n$ be a polynomial and $I \subset \aru n$ an ideal.
We write $f_{h_1}=\Phi_{h_1}(f) \in \aru {n-1}$ and $I_{h_1}$ for the ideal $\Phi_{h_1}(I) =\{ \Phi_{h_1}(f): f \in I\}$ of $\aru {n-1}$.
Let $h_2$ be an another linear form of $\aru n$.
Assume that the coefficient of $x_{n-1}$ in $(h_2)_{h_1}$ is not zero.
Then define $f_{\langle h_1,h_2\rangle } = \Phi_{(h_2)_{h_1}}(f_{h_1})$ and $I_{\langle h_1,h_2\rangle} = \Phi_{(h_2)_{h_1}} (I_{h_1})$.
Inductively, we define $I_{\langle h_1,h_2,\dots,h_m\rangle}$ by the same way
for linearly independent linear forms $h_1,h_2,\dots,h_m$ of $\aru n$,
where we assume that the coefficient of $x_{n+1-k}$ in $(h_k)_{\langle h_1,\dots,h_{k-1}\rangle}$ is not zero for each $1 \leq k\leq m$.

\begin{lemma}[{\cite[Corollary 2.15]{G}}]  \label{hyper}
Let $I \subset \aru n$ be a graded ideal and  $h_1,\dots,h_m$ generic linear forms of $\aru n$ with $1\leq m \leq n$.
Then
$$\gin(I_{\langle h_1,\dots,h_m \rangle })=\gin(I)_{\langle x_{n-m+1},\dots,x_n \rangle} = \gin(I) \cap \aru {n-m}.$$
\end{lemma}

\begin{theorem} \label{sqinc}
Let $K$ be a field of characteristic $0$ and $\Gamma$ a $(d-1)$-dimensional Gorenstein* complex on $[n]$
with the weak Lefschetz property.
Then, one has
$$\beta_{ij}^{\aru n} (I_\Gamma) \leq \beta_{ij}^{\aru n} (I_{{\mathrm{Sq}(\Gamma)}}) \hspace{10pt} \mbox{ for all } i,j.$$
\end{theorem}

\begin{proof}
Let $\vartheta_1,\vartheta_2,\dots,\vartheta_d,\omega$ be generic linear forms of $\aru n$.
Then Lemma \ref{glf} says that
$\vartheta_1,\vartheta_2,\dots,\vartheta_d$ is a system of parameters of $\aru n / I_\Gamma$ 
and $\omega$ is a weak Lefschetz element of $(\aru n /I_\Gamma) \otimes ( \aru n/\langle \vartheta_1,\dots,\vartheta_d \rangle)$.

Let $\tilde R = (\aru n /\langle \vartheta_1,\vartheta_2,\dots,\vartheta_d,\omega\rangle)$
and $A= (\aru n/ I_\Gamma) \otimes (\aru n/\langle\vartheta_1,\vartheta_2,\dots,\vartheta_d\rangle).$
Then, by the definition of $(I_\Gamma)_{\langle \vartheta_1,\dots,\vartheta_d,w\rangle}$, we have 
\begin{eqnarray}
\beta_{ij}^{\tilde R}(A / \omega A) = \beta_{ij}^{\tilde R}( \tilde R/ (I_\Gamma \otimes \tilde R))
= \beta_{ij}^{\aru {n-d-1}} ( \aru {n-d-1}/( (I_\Gamma)_{\langle \vartheta_1,\dots,\vartheta_d,w\rangle})) \label{sono1} 
\end{eqnarray}
for all $i,j$.
Recall that $I(U(\Gamma))$ is the ideal of $\aru n$ generated by all monomials $u \in \aru {n-d-1}$
with $u \in \gin(I_\Gamma)$.
Thus $\gin(I_\Gamma) \cap \aru {n-d-1} = I(U(\Gamma)) \cap \aru {n-d-1}$.
Also, it is known that $\beta_{ij}^{\aru n} (\aru n/ I) \leq \beta_{ij}^{\aru n} (\aru n / \gin (I))$ for any graded ideal $I \subset \aru n$
(see e.g., \cite[Theorem 3.1]{H}).
Thus, by Lemma \ref{hyper},
we have
\begin{eqnarray*}
\beta_{ij}^{\aru {n-d-1}} ( \aru {n-d-1}/( (I_\Gamma)_{\langle \vartheta_1,\dots,\vartheta_d,w\rangle})) &\leq&
\beta_{ij}^{\aru {n-d-1}} ( \aru {n-d-1}/ \gin ((I_\Gamma)_{\langle \vartheta_1,\dots,\vartheta_d,w\rangle}))\\ 
&=& \beta_{ij}^{\aru {n-d-1}} ( \aru {n-d-1}/( \gin (I_\Gamma )\cap \aru {n-d-1})) \\ 
&=& \beta_{ij}^{\aru {n-d-1}} ( \aru {n-d-1}/ (I(U(\Gamma ))\cap \aru {n-d-1})) 
\end{eqnarray*}
for all $i,j$.
By the definition of $I(U)$, we have $\max(I(U)) \leq n-d-1$.
Thus we have $$\beta_{ij}^{\aru {n-d-1}} ( \aru {n-d-1}/ (I(U(\Gamma ))\cap \aru {n-d-1}))=
 \beta_{ij}^{\aru n} (\aru n/I(U(\Gamma)))$$
 for all $i,j$.
Then the equality (\ref{sono1}) together with the above computations says
$$\beta_{ij}^{\tilde R} (A / \omega A) \leq \beta_{ij}^{\aru {n-d-1}} (\aru {n-d-1} / (I(U(\Gamma)) \cap \aru {n-d-1}))
= \beta_{ij}^{\aru n} (\aru n/I(U(\Gamma)))$$ 
for all $i,j$.
Then, by Lemma \ref{gbetti}, we have
$$\beta_{ii+j}^{\aru n}( \aru n /I_\Gamma)
\leq \beta_{ii+j}^{\aru n} (\aru n/I(U(\Gamma)))
+ \beta_{n-d-i,n-i-j}^{\aru n} (\aru n/I(U(\Gamma)))$$
for all $i,j$.
Since $I(U(\Gamma))$ has no generators of degree $j > \lfloor \frac d 2 \rfloor$ +1,
it follows from the Eliahou--Kervaire formula that
$\beta_{ii+j}^{\aru n} (\aru n/I(U(\Gamma)))=0$ for $j > \frac d 2$.
Thus we have
\begin{eqnarray*}
\beta_{ii+j}^{\aru n}( \aru n / I_\Gamma) \leq
\left\{
\begin{array}{ll}
\beta_{ii+j}^{\aru n}(\aru n / I(U(\Gamma))), \hspace{138pt}
& \mbox{for} \ j  <  \frac {d} 2, \\
\beta_{ii+j}^{\aru n}(\aru n/ I(U(\Gamma))) 
+\beta_{n-d-i,n-i-j}^{\aru n}(\aru n/ I(U(\Gamma))), \hspace{10pt} 
&\mbox{for } j= \frac d 2, \\
\beta_{n-d-i,n-i-j}^{\aru n}(\aru n / I(U(\Gamma))), \hspace{100pt}
&\mbox{for } j >  \frac d 2.
\end{array}
\right.
\end{eqnarray*}
Since $\mathrm{Sq(\Gamma)}=S_d(U(\Gamma))$,
Theorem \ref{sqbetti} says that
$\beta_{ij}^{\aru n}(\aru n /I_\Gamma) \leq \beta_{ij}^{\aru n} (\aru n/I_{\mathrm{Sq}(\Gamma)})$ for all $i$ and $j$.
\end{proof}

\begin{exam} 
Let $U=\{1,x_1,x_2,x_3,x_1x_3,x_2x_3,x_3^2\}$.
Then the squeezed $5$-ball $B_5(U)$ is the simplicial complex on $\{1,2,\dots,9\}$ generated by
\begin{eqnarray*}
F_5(U)=\left\{
\begin{array}{l}
\{4,5,6,7,8,9\},\{1,2,6,7,8,9\},\{2,3,6,7,8,9\},\{3,4,6,7,8,9\},\\
\{1,2,5,6,8,9\},\{2,3,5,6,8,9\},\{3,4,5,6,8,9\}
\end{array}
\right\}
\end{eqnarray*}
and $I(U)=\langle x_1^2,x_1x_2,x_2^2,x_1x_3^2,x_2x_3^2,x_3^3 \rangle$.
Let $R=K[x_1,x_2,\dots,x_9]$.
Then Proposition \ref{ginball} guarantees that
$$I_{B_5(U)}=\langle x_1x_3,x_1x_4,x_2x_4,x_1x_5x_7,x_2x_5x_7,x_3x_5x_7\rangle $$
and the minimal graded free resolution of $I_{B_5(U)}$ is of the form
\[
0 \longrightarrow 
R(-5)^3 \longrightarrow 
R(-4)^6 \oplus R(-3)^2 \longrightarrow 
R(-3)^3\oplus R(-2)^3 \longrightarrow 
I_{B_5(U)} \longrightarrow  0.
\]
Also, Theorem \ref{sqbetti} says that the minimal graded free resolution of $R/I_{S_5(U)}$ is of the form
\begin{eqnarray*}
0 \to  R(-9)
 &\longrightarrow&
R(-7)^3 \oplus R(-6)^3 \oplus R(-5)^3  \\ 
&\longrightarrow& 
R(-6)^2 \oplus R(-5)^6 \oplus R(-4)^6 \oplus R(-3)^2 \\ 
&\longrightarrow& 
\hspace{53pt}R(-4)^3 \oplus R(-3)^3 \oplus R(-2) ^3
\longrightarrow 
R
\longrightarrow
R/I_{S_5(U)}
\to 0.
\end{eqnarray*}
Note that
$I_{S_5(U)}=I_{B_5(U)} +\langle x_2x_6x_8x_9,x_3x_6x_8x_9,x_4x_6x_8x_9 \rangle.$
\end{exam}

\begin{exam}\label{referee}
We will give an easy example of a simplicial sphere whose graded Betti numbers strictly increase
by squeezing.

Let $\Gamma$ be the boundary complex of the octahedron.
Then $n=6$, $d=3$, $I_\Gamma = (x_1x_2,x_3x_4,x_5x_6)$ and
the minimal graded free resolution of $R/I_\Gamma$ is of the form
\begin{eqnarray*}
0 \longrightarrow  R(-6)
\longrightarrow
R(-4)^3  
\longrightarrow
R(-2)^3
\longrightarrow
R
\longrightarrow
R/I_\Gamma
\longrightarrow
 0.
\end{eqnarray*}
On the other hand, \vspace{1pt}
since $\lfloor \frac d 2 \rfloor =1$ and $n-d-1=2$, we have $U(\Gamma)=\{1,x_1,x_2\}$
and $I(U(\Gamma))=(x_1^2,x_1x_2,x_2^2)$.
Since $\mathrm{Sq}(\Gamma)=S_d(U(\Gamma))$, the minimal graded free resolution of $R/I_{\mathrm{Sq(\Gamma)}}$ is of the form
\begin{eqnarray*}
0 \longrightarrow  R(-6)
&\longrightarrow&
R(-4)^3 \bigoplus R(-3)^2 \\ 
&\longrightarrow&
R(-3)^2 \bigoplus R(-2)^3
\longrightarrow
R
\longrightarrow
R/I_{\mathrm{Sq(\Gamma)}}
\longrightarrow
0.
\end{eqnarray*}

\end{exam}

\section{Characterization of generic initial ideals associate with simplicial $d$-polytopes for $d \leq 5$}

We refer the reader to \cite{Gr} for the foundations of convex polytopes.
A $(d-1)$-dimensional simplicial complex $\Gamma$ is called \textit{polytopal} if 
$\Gamma$ is isomorphic to the boundary complex of a simplicial $d$-polytope.
Although most of the squeezed $(d-1)$-spheres are not polytopal for $d \geq 5$,
Pfeifle proved that squeezed $3$-spheres are polytopal.
His proof implies the following fact.

\begin{lemma} \label{4polytopal}
S-squeezed $4$-spheres are polytopal.
\end{lemma}

\begin{proof}
We recall Pfeifle's proof
(see  \cite[pp. 400--401]{P}).
Let $\mathcal{C}_5(n-1)$ be the collection of facets of the boundary complex of the cyclic $5$-polytope with $n-1$ vertices.
Let $U\subset \aru {n-6}$ be a shifted order ideal of monomials of degree at most $2$.
Recall that $F_4(U)=\{ F_4(u)\subset [n-1]:u \in U\}$ can be regarded as a subcollection of $\mathcal{C}_5(n-1)$.
We identify each $F \in \mathcal{C}_5(n-1)$ and the corresponding facet of the cyclic $5$-polytope with $n-1$ vertices.

Pfeifle proved that there is a set $C=\{v_1,v_2,\dots,v_{n-1}\}$ of vertices on $\mathbb{R}^5$
and there is a vertex $v_n$ on $\mathbb{R}^5$
such that $\mathrm{conv}(C)$ is the cyclic $5$-polytope with $(n-1)$ vertices and
\begin{itemize}
\item[(i)]
$v_n$ is beyond $F$ for $F \in F_4(U)$;
\item[(ii)]
$v_n$ is beneath $F$ for $F \in \mathcal{C}_5(n-1) \setminus F_4 (U)$.
\end{itemize} 
See \cite[\S 5.2]{Gr} for the definitions of beneath and beyond.

Let $H$ be the hyperplane separating $v_n$ from $C$
with equation $\langle a,x \rangle=a_0$,
where $a\in \mathbb{R}^5$, $a_0 \in \mathbb{R}$ and $\langle u,v \rangle$ is the scalar product of $u,v \in \mathbb{R}^5$.
Assume that $v_n$ is sufficiently close to $H$.
Consider the projective transformation $\varphi$ defined by
$\varphi (x)= \frac {x} {( \langle a,x \rangle -a_0)}$.
Then $\mathrm{conv}(\varphi(C))$ is isomorphic to $\mathrm{conv} (C)$ since all vertices in $C$
lies on the same side of $H$.
However, by the projective transformation $\varphi$,
when we regard $\mathcal{C}_5(n-1)$ as the collection of facets of $\mathrm{conv}(\varphi(C))$,
the vertex $\varphi(v_n)$ becomes  beneath $F$ for $F \in F_4(U)$
and beyond $F$ for $F \in \mathcal{C}_5(n-1) \setminus F_4 (U)$.

We will explain why this occurs.
Let $H(t)$ be the hyperplane with equation $\langle a,x \rangle =a_0+t$.
Set $M= \frac {1} {\langle a,v_n \rangle -a_0}$ and assume $M>0$.
Since $v_n$ is sufficiently close to $H$,
there exist $\delta\in \mathbb{R}$ such that
$|\delta| \ll \min_{v_k \in C} \{ |\langle a,v_k \rangle -a_0|\}$,
$\frac{1}{|\langle a,v_n \rangle -a_0-\delta|} \geq M$
and $\langle a,v_n\rangle-a_0-\delta<0$.
Then $C$ and $v_n$ lie on the same side of $H(\delta)$.
Let $\varphi'$ be the projective transformation defined by
$\varphi'(x)=\frac {x} {\langle a,x\rangle -a_0-\delta}.$
Then $\mathrm{conv}( \varphi'(C) \cup \{\varphi'(v_n)\})$ is isomorphic to $\mathrm{conv}(C \cup \{v_n\})$
and $\mathrm{conv}( \varphi'(C))$ is isomorphic to $\mathrm{conv}( \varphi(C))$.
On the other hand,
since 
$|\delta|$ is sufficiently small,
the difference between $\mathrm{conv}( \varphi'(C)) $ and $\mathrm{conv}(\varphi(C))$ is also sufficiently small.
Let $M'= \frac {1} {\langle a,v_n\rangle -a_0-\delta}$.
Then we have $\varphi(v_n)=M v_n$, $\varphi'(v_n)=M' v_n$ and $|M| \leq |M'|$.
Since $M>0$, $M'<0$ and $|M|$ is sufficiently large,
if $\varphi'(v_n)$ is beneath (beyond) for a facet $F$ of $\mathrm{conv}( \varphi'(C))$,
then $\varphi(v_n)$ is beyond (beneath) $F$ when we regard $F$ as a facet of $\mathrm{conv}( \varphi(C))$,
as required.
\medskip

Let $P= \mathrm{conv}( \varphi(C) \cup \{ \varphi (v_n)\})$.
We will show $S_5(U)$ is the boundary complex of the simplicial $5$-polytope $P$.
Since $B_5(U)$ is generated by $\{ \{n\} \cup F_4(u):u\in U\}$,
it follows from \cite[Proposition 1]{L} that $S_5(U)$ is the simplicial complex generated by
$$F_4(U) \cup \{ F \cup \{n\} : F \mbox{ is a facet of } S_4(U)\}.$$
Recall that $B_4(U)$ is the $5$-ball generated by $F_4(U)$ and  $S_4(U)$ is its boundary.
It follows form \cite[\S 5.2 Theorem 1]{Gr} that
$F$ is a facet of the boundary complex of $P$ with $\{n\} \not \in F$ if and only if 
$F$ is a facet of the boundary complex of $\mathrm{conv}(\varphi(C))$ and $\varphi(v_n)$ is beneath $F$.
Also, since $S_4(U)$ is the boundary of the $5$-ball generated by $F_4(U)$,
it follows that $F$ is a facet of the boundary complex of $P$ with $\{n \} \in F$ if and only if
$F \setminus \{n\}$ is a facet of $S_4(U)$.
Thus $S_5(U)$ is the boundary complex of the simplicial $5$-polytope $P$.
\end{proof}

If $\Gamma$ is the boundary complex of a simplicial $d$-polytope on $[n]$,
then $\Gamma$ is a $(d-1)$-dimensional Gorenstein* complex with the strong Lefschetz property
(see \cite[pp. 75--78]{St}).
Thus, as we saw at the end of \S 4,
Lemma \ref{4polytopal} yields a complete characterization of generic initial ideals of Stanley--Reisner ideals
of the boundary complexes of simplicial $d$-polytopes 
(or $(d-1)$-dimensional Gorenstein* complex with the strong Lefschetz property)
for $d \leq 5$ when the base field is of characteristic $0$.

\begin{theorem} \label{chara5}
Let $K$ be a field of characteristic $0$ and $d \leq 5$.
Let $I\subset \aru n$ be a strongly stable ideal and
$A=\aru{n-d} /(I \cap \aru{n-d})$.
Then there is the boundary complex $\Gamma$ of a simplicial $d$-polytope on $[n]$
such that $I = \gin(I_\Gamma)$ if and only if 
$\max(I) = {n-d}$, $A_{d+1} =\{0\}$, $\dim_K A_1=n-d$ and the multiplication map
$x_{n-d}^{d-2i}:A_i \to A_{d-i}$ is an isomorphism  for $0\leq i \leq \lfloor \frac d 2 \rfloor$.
\end{theorem}

\begin{proof}
Let $A=\aru {n-d} / (\gin(I_\Gamma) \cap \aru {n-d})$.
Then $\dim_K A_1=n-d$ is obvious.
Since $\Gamma$ is a $(d-1)$-dimensional Gorenstein* complex with the strong Lefschetz property,
Lemmas \ref{hcm} and \ref{lefschetz1} say that
$\max(I_\Gamma) = {n-d}$, $A_{d+1} = \{0\}$, $\dim_K A_1=n-d$ and the multiplication map
$x_{n-d}^{d-2i}:A_i \to A_{d-i}$ is an isomorphism  for $0\leq i \leq \lfloor \frac d 2 \rfloor$.

Conversely,
given a strongly stable ideal $I\subset \aru n$ which satisfies the conditions of Theorem \ref{chara5}.
Then,
$U=\{ u\in \aru {n-d-1}: u\not\in I \mbox{ is a monomial}\}$
is a shifted order ideal of monomials of degree at most $\lfloor \frac d 2 \rfloor$
and $U$ determines $I$ in the same way as Lemma \ref{lefschetz2}.
Then what we must do is finding the boundary complex $\Gamma$ of a simplicial $d$-polytope with $U(\Gamma) =U$.
Now, Theorem \ref{conj} says $U(S_d(U))=U$.
Since $d\leq 5$, this $S_d(U)$ is polytopal by Lemma \ref{4polytopal}.
\end{proof} 

Theorem \ref{chara5} is not true for $d \geq 6$.
Let $sq(d,n)$ be the number of squeezed $(d-1)$-spheres on $[n]$,
$ssq(d,n)$ the number of S-squeezed $(d-1)$-spheres on $[n]$
and $c(d,n)$ the number of combinatorial type of the boundary complex of simplicial $d$-polytopes with $n$ vertices.
Then it is known that 
$$ \log(c(d,n)) \leq d(d+1)n \log (n) $$
and 
$$ \log(sq (d,n)) \geq \frac {1} {(n-d)(d+1)} {{ n+ \lfloor \frac {d+2} 2 \rfloor} \choose {\lfloor \frac {d+1} 2 \rfloor}}.$$
( See \cite{K3} or \cite[pp. 397]{P}.)
Thus we have $c(d,n) \ll sq (d,n)$ for $d \geq 5$ and for $n \gg 0$.
On the other hand,
it is easy to see that 
$sq(d-1,n-1) =ssq(d,n)$ and 
$ssq(d,n)$ is equal to the number of shifted order ideals $U \subset \aru {n-d-1}$
of monomials of degree at most $\lfloor \frac {d} 2 \rfloor$.
Then, the above upper bound for $c(d,n)$ and the lower bound for $sq(d,n)$ imply
$c(d,n) \ll ssq (d,n)$ for $d \geq 6$ and for $n \gg 0$.
Thus the number of strongly stable ideals which satisfies the condition of Theorem \ref{chara5} is strictly larger than
the number of combinatorial type of the boundary complex of simplicial $d$-polytopes with $n$ vertices
for $d \geq 6$ and for $n \gg 0$.

\section{Exterior algebraic shifting of squeezed balls}

Let $K$ be an infinite field,
$V$ a $K$-vector space of dimension $n$ with basis $e_1,\dots,e_n$
and $E=\bigoplus_{d=0}^n \bigwedge^d V$ the exterior algebra of $V$.
For a subset $S = \{ s_1, s_2,\dots,s_k \} \subset [n]$ with $s_1\leq s_2\leq \cdots \leq s_k$,
we write $e_S=e_{s_1}\wedge e_{s_2} \wedge \cdots \wedge e_{s_k}\in E$ and $x_S=x_{s_1}x_{s_2} \cdots x_{s_k} \in \aru n$.
The element $e_S$ is called the \textit{monomial} of $E$ of degree $k$.
In the exterior algebra,
the generic initial ideal $\Gin(J)$ of a graded ideal $J\subset E$
is defined similarly as in the case of the polynomial ring 
(\cite[Theorem 1.6]{AHH}).

Let $\Gamma$ be a simplicial complex on $[n]$.
The \textit{exterior face ideal} $J_\Gamma\subset E$ of $\Gamma$  is the monomial ideal 
generated by all monomials $e_S\in E$ with $S \not\in \Gamma$.
The \textit{exterior algebraic shifted complex} $\dele \Gamma$ of $\Gamma$ is the simplicial complex on $[n]$ defined by
$J_{\dele \Gamma}=\Gin(J_\Gamma)$.
Thus knowing $\dele \Gamma$ is equivalent to knowing $\Gin(J_\Gamma)$.

A squarefree monomial ideal $I \subset \aru n$ is called \textit{squarefree strongly stable}
if $v \prec u$ and $u \in I$ imply $v \in I$ for all squarefree monomials $u$ and $v$ in $\aru n$.
A simplicial complex $\Gamma$ on $[n]$ is called \textit{shifted} if $I_\Gamma\subset \aru n$ is squarefree strongly stable.
We recall basic properties of $\Delta^e$.

\begin{lemma}[{\cite[Proposition 8.8]{H}}] \label{itumono}
Let $\Gamma$ and $\Gamma '$ be simplicial complexes on $[n]$.
Then
\begin{itemize}
\item[(i)] $I_{\dele \Gamma}$ is squarefree strongly stable;
\item[(ii)]
$I_\Gamma$ and $I_{\dele \Gamma}$ have the same Hilbert function;
\item[(iii)]
if $I_\Gamma \subset I_\Gamma'$, then $I_{\dele \Gamma} \subset I_{\dele{\Gamma'}}$.
\end{itemize}
\end{lemma}

Let $\Gamma$ be a simplicial complex on $[n]$.
The \textit{cone} $Cone(\Gamma,n+1)$ over $\Gamma$ is the simplicial complex on $[n+1]$
generated by $\{ \{n+1\}\cup S: S \in \Gamma\}$.
In other words, $Cone(\Gamma,n+1)$ is the simplicial complex defined by
$I_{Cone(\Gamma,n+1)}= I_\Gamma \aru {n+1}$.

\begin{lemma}[{\cite[Corollary 5.5]{N}}] \label{cone}
Let $n > m > 0$ be positive integers,
$\Gamma$ a simplicial on $[n]$ and $\Gamma'$ a simplicial complex on $[m]$.
If $I_\Gamma = I_{\Gamma'} \aru n$, then $I_{\dele \Gamma} = I_{\dele {\Gamma'}}\aru n$.
\end{lemma}

\begin{proof}
If $I_\Gamma = I_{\Gamma'} \aru n$, then $\Gamma$ is obtained from $\Gamma'$
by taking a cone repeatedly.
On the other hand,
it follows from \cite[Corollary 5.5]{N} that
$$\dele {Cone (\Gamma' ,n+1)} = Cone ( \dele {\Gamma'},n+1).$$
Thus the assertion follows.
\end{proof}


We will prove the analogue of Theorem \ref{soperator} for generic initial ideals in the exterior algebra.
A map $\sigma: \setmono \to \setmono$ is called \textit{a squarefree stable operator}
if $\sigma$ satisfies
\begin{itemize}
\item[(i)]
if $I\subset \aru \infty$ is a finitely generated squarefree strongly stable ideal,
then $\sigma(I)$ is also a squarefree monomial ideal and $\beta_{ij}(I) = \beta_{ij}(\sigma(I))$ for all $i,j$;
\item[(ii)]
if $J \subset I$ are finitely generated squarefree strongly stable ideals of $\aru \infty$,
then $\sigma (J) \subset \sigma(I)$.
\end{itemize}

Like strongly stable ideals, the graded Betti numbers of a squarefree strongly stable ideal $I \subset \aru n$ are
given by the formula (\cite[Corollary 3.6]{H})
\begin{eqnarray}
\beta_{i,i+j}(I) = \sum_{u\in G(I), \ \deg(u)=j} { m(u) - j \choose i}. \label{sqformula}
\end{eqnarray}

Let $\sigma: \setmono \to \setmono$ be a squarefree stable operator and $\Gamma$ a simplicial complex on $[n]$.
We write $\sigma (\Gamma)$ for the simplicial complex on $[n]$ with
$$I_{\sigma(\Gamma)}= \sigma(I_\Gamma \aru \infty) \cap \aru n.$$

\begin{lemma} \label{hodai}
Let $\sigma: \setmono \to \setmono$ be a squarefree stable operator,
$\Gamma$ a shifted simplicial complex on $[n]$.
Assume $n \geq \max(\sigma(I_\Gamma\aru \infty))$.
Then one has 
$\max(I_\Gamma) = \max(I_{\dele {\sigma(\Gamma)}})$.
In particular, for all $ n \geq m \geq \max(I_\Gamma)$ and for all $d\geq 0$, one has
$$|\{x_S\in (I_\Gamma)_d: S \subset [m]\}|=|\{x_S \in (I_{\dele {\sigma(\Gamma)}})_d: S \subset [m]\}|.$$
\end{lemma}

\begin{proof}
First, we will show $\max(I_{\dele {\sigma(\Gamma)}}) =\max(I_\Gamma)$.
The formula (\ref{sqformula}) says that,
for every squarefree strongly stable ideal $J\subset \aru n$,
one has
\begin{eqnarray}
\max(J)= \max\{k: \beta_{ik}(J)\ne 0 \mbox{ for some } i\}. \label{hosi71}
\end{eqnarray}
Also, it follows from \cite[Theorem 7.1]{H} that
\begin{eqnarray*}
\max\{k: \beta_{ik}(I_{\sigma(\Gamma)}) \ne 0\mbox{ for some } i\}
=\max\{k: \beta_{ik}(I_{\dele{\sigma(\Gamma)}})\ne 0\mbox{ for some } i\}.
\end{eqnarray*}
Since $n \geq \max(\sigma (I_\Gamma \aru \infty))$,
$I_\Gamma$ and $I_{\sigma(\Gamma)}= \sigma (I_\Gamma \aru \infty) \cap \aru n$ have the same graded Betti numbers.
Thus we have
\begin{eqnarray}
\max\{k: \beta_{ik}(I_\Gamma)\ne 0\mbox{ for some } i\}
=\max\{k: \beta_{ik}(I_{\dele{\sigma(\Gamma)}})\ne 0\mbox{ for some } i\}. \label{hosi72}
\end{eqnarray}
Since $I_\Gamma$ and $I_{\dele {\sigma (\Gamma)}}$ are squarefree strongly stable,
the equalities (\ref{hosi71}) and (\ref{hosi72}) say $\max(I_{\dele{\sigma(\Gamma)}})=\max(I_\Gamma)$.

Since $I_\Gamma$ and $I_{\sigma(\Gamma)}$ have the same graded Betti numbers,
Lemma \ref{itumono} says that
$I_\Gamma$, $I_{\sigma(\Gamma)}$ and $I_{\dele {\sigma(\Gamma)}}$ have the same Hilbert function.
Thus $I_\Gamma \cap \aru m$ and $I_{\dele{\sigma(\Gamma)}} \cap \aru m$ have the same Hilbert function
for all $ n \geq m \geq \max(I)$.
Since the Hilbert function of Stanley--Reisner ideal $I_\Gamma$ of $\Gamma$ is determined by the $f$-vector of $\Gamma$,
the previous fact says 
$$|\{x_S\in (I_\Gamma)_d: S \subset [m]\}|=|\{x_S \in (I_{\dele {\sigma(\Gamma)}})_d: S \subset [m]\}|$$
for all $d \geq 0$ and for all $ n \geq m \geq \max(I)$.
\end{proof}

\begin{proposition} \label{sqoperator}
With the same notation as in Lemma \ref{hodai}.
One has $I_{\dele {\sigma(\Gamma)}} = I_\Gamma$.
\end{proposition}

\begin{proof}
Let $m =\max(I_\Gamma)$.
By virtue of Lemma \ref{hodai}, what we must prove is
$$I_\Gamma \cap \aru m = I_{\dele {\sigma (\Gamma)}} \cap \aru m.$$

We use induction on $m$.
In case of $m=1$, since $G(I_\Gamma)=\{x_1\}$
and  since $I_\Gamma \cap K[x_1]$ and $I_{\dele {\sigma (\Gamma)}} \cap K[x_1]$ have the same Hilbert function,
we have $G( I_{\dele {\sigma (\Gamma)}})=\{x_1\}$.

Assume $m >1$.
Fix an integer $d\geq 0$.
Write $I_{\langle d \rangle}\subset \aru \infty$ for the ideal generated by all squarefree monomials $u \in I_\Gamma \cap \aru m$
of degree $d$. 
Consider the colon ideal $J= ( I_{\langle d \rangle }:x_m)\subset \aru \infty$.
Then $I_{\langle d \rangle}$ and $J$ are also squarefree strongly stable and $\max(J)<m$.

Let $l= \max \{ \max (\sigma(J)),\max(\sigma(I_{\langle d \rangle})),n \}$.
Let $\Gamma'$ and $\Gamma''$ be simplicial complexes on $[l]$
with $I_{\Gamma'}=I_{\langle d \rangle }\cap \aru l$ and with  $I_{\Gamma''}=J\cap \aru l$.
Since $\max(J)<m$,
the assumption of induction says $I_{\Gamma''}= I_{\dele {\sigma(\Gamma'')}}$.
Also, since $I_{\Gamma''} \supset I_{\Gamma'}$ are squarefree strongly stable ideals,
we have $I_{\sigma(\Gamma'')} \supset I_{\sigma(\Gamma')}$.
Thus Lemma \ref{itumono} (iii) says
\begin{eqnarray}
I_{\Gamma ''}
=I_{\dele {\sigma (\Gamma'')}}  \supset I_{\dele {\sigma (\Gamma ')}}. \label{hosihosi3}
\end{eqnarray}
Let $\Sigma$ be the simplicial complex on $[l]$ with $I_\Sigma = I_\Gamma \aru l$.
Since $l \geq n$, we have $I_{\sigma(\Sigma)} = I_{\sigma( \Gamma)} \aru l.$
Thus, by Lemma \ref{cone}, we have
$$I_{\dele {\sigma (\Sigma)}} = I_{\dele {\sigma (\Gamma)}} \aru l.$$
Also, since $I_\Gamma \aru \infty \supset I_{\langle d \rangle}$,
we have $I_\Sigma =(I_\Gamma \aru \infty) \cap \aru l \supset I_{\langle d \rangle }\cap \aru l= I_{\Gamma'}$.
Note that $I_\Sigma$ is squarefree strongly stable.
Then we have
\begin{eqnarray*}
I_{\dele  {\sigma(\Gamma)}} \aru l =I_{\dele {\sigma (\Sigma)}} \supset I_{\dele {\sigma (\Gamma')}}.
\end{eqnarray*}
In particular,
since $m \leq n \leq l$, we have
\begin{eqnarray}
I_{\dele  {\sigma(\Gamma)}} \cap \aru m  \supset I_{\dele {\sigma (\Gamma')}}\cap \aru m. \label{hosihosi4}
\end{eqnarray}

Next, we will show
\begin{eqnarray}
\{ x_S \in (I_{\Gamma''})_d: S \subset [m]\}
=\{ x_S \in (I_{\Gamma'})_d: S \subset [m]\}
=\{ x_S \in (I_{\Gamma})_d: S \subset [m]\}. \label{hosihosi1}
\end{eqnarray}
The second equality directly follows from the definition of $\Gamma'$.
Also, $I_{\Gamma''} \supset I_{\Gamma'}$ is obvious.
Thus what we must prove is 
$\{ x_S \in (I_{\Gamma''})_d: S \subset [m]\}\subset \{ x_S \in (I_{\Gamma'})_d: S \subset [m]\}.$
Let $x_S \in I_{\Gamma''} \cap \aru m$ be a squarefree monomial.
Then $x_Sx_m \in I_{\langle d\rangle }\cap \aru m$.
Since $I_{\langle d \rangle}$ is squarefree monomial ideal,
there is a squarefree monomial $x_T \in I_{\langle d \rangle } \cap \aru m$ of degree $d$
such that $x_T$ divides $x_S x_m$.
Then we have $x_T = (x_S x_m) / x_i $ for some $i \in S \cup \{m\}$.
Thus we have $x_S \preceq x_T$.
Since $I_{\langle d\rangle } $ is squarefree strongly stable,
we have $x_S \in I_{\langle d \rangle } \cap \aru m =I_{\Gamma'} \cap \aru m$.

Lemma \ref{hodai} together with (\ref{hosihosi1}) says 
\begin{eqnarray*}
|\{ x_S \in (I_{\dele {\sigma(\Gamma'')}})_d: S \subset [m]\}|
&=&|\{ x_S \in (I_{\dele {\sigma(\Gamma')}})_d: S \subset [m]\}|\\
&=&|\{ x_S \in (I_{\dele {\sigma(\Gamma)}})_d: S \subset [m]\}|. 
\end{eqnarray*}
Then, the equalities (\ref{hosihosi3}) and  (\ref{hosihosi4}) together with (\ref{hosihosi1}) says
\begin{eqnarray*}
\{ x_S \in (I_{\Gamma})_d: S \subset [m]\}=
\{ x_S \in (I_{\Gamma''})_d: S \subset [m]\}
&=&\{ x_S \in (I_{\dele {\sigma(\Gamma')}})_d: S \subset [m]\}\\
&=&\{ x_S \in (I_{\dele {\sigma(\Gamma)}})_d: S \subset [m]\}.
\end{eqnarray*}
Thus, for any squarefree monomial $x_S \in \aru m$,
we have $x_S \in I_\Gamma \cap \aru m$ if and only if 
$x_S \in I_{\dele {\sigma (\Gamma)}} \cap \aru m$.
Hence we have $I_\Gamma\cap \aru m = I_{\dele {\sigma (\Gamma)}} \cap \aru m$ as required.
\end{proof}

Next, we will show that the maps $\alpha^a:\setmono \to \setmono$, which we define in \S 1,
are squarefree stable operators.

\begin{lemma}[{\cite[Lemmas 8.17 and 8.20]{H}}] \label{sss}
Let $\alpha:\setmono \to \setmono$ be the map defined by
$$\alpha(x_{i_1}x_{i_2} \cdots x_{i_k})=x_{i_1} x_{i_2+1} \cdots x_{i_k +k-1}$$
for any monomial $x_{i_1}x_{i_2} \cdots x_{i_k} \in \setmono$ with $i_1 \leq i_2 \leq \dots \leq i_k$.
\begin{itemize}
\item[(i)]
If $I \subset \aru \infty$ is a finitely generated squarefree strongly stable ideal,
then there is the strongly stable ideal $I' \subset \aru \infty$ such that $\alpha(I')=I$.
\item[(ii)]
If $I\subset \aru \infty$ is a finitely generated strongly stable ideal,
then $\alpha(I)\subset \aru \infty$ is a squarefree strongly stable ideal.
\end{itemize}
\end{lemma}

\begin{proposition} \label{sqsatisfy}
Let $a=(0,a_1,a_2,a_3,\dots)$ be a nondecreasing infinite sequence of integers.
Let $\alpha^a: \setmono \to \setmono$ be the map defined by
$$\alpha^a(x_{i_1}x_{i_2}x_{i_3} \cdots x_{i_k}) =x_{i_1} x_{i_2+a_1} x_{i_3 +a_2}\cdots x_{i_k +a_{k-1}}$$
for any monomial $x_{i_1}x_{i_2} \cdots x_{i_k}\in \setmono$ with $i_1 \leq i_2 \leq \dots \leq i_k$.
Then $\alpha^a :\setmono \to \setmono $ is a squarefree stable operator.
\end{proposition}

\begin{proof}
First, we will prove the property (i) of squarefree stable operators.
Let $\alpha: \setmono \to \setmono$ be the map in Lemma \ref{sss} and
let $a'=(0,a_1+1,a_2+2,a_3+3,\dots)$.
For any finitely generated squarefree strongly stable ideal $I \subset \aru \infty$,
Lemma \ref{sss} says that there is the strongly stable ideal $I' \subset \aru \infty$
such that $\alpha(I') =I$.
Since
$\alpha^a(I)=\alpha^a(\alpha(I'))=\alpha^{a'}(I')$
and since Proposition \ref{satisfy} says that $\alpha$ and $\alpha^{a'}$ are stable operators,
we have
$$\beta_{ij} (I)=\beta_{ij}(I')= \beta_{ij}(\alpha^{a'}(I')) =\beta_{ij} (\alpha^a(I))$$
for all $i,j$.

Second, we will prove the property (ii) of squarefree stable operators.
If $I \supset J$ are finitely generated squarefree strongly stable ideals of $\aru \infty$,
then, for any $u=x_{i_1}x_{i_2} \cdots x_{i_k} \in G(J)$ with $i_1 < i_2 < \dots < i_k$,
there is $w \in G(I)$ such that $w$ divides $u$.
Since $I$ is squarefree strongly stable,
we may assume $w = x_{i_1}x_{i_2}\cdots x_{i_l}$ for some $l \leq k$.
Then $\alpha^a(w)$ divides $\alpha^a(u)$.
Thus we have $\alpha^a(u) \in \alpha^a(I)$ and, therefore, $\alpha^a(I) \supset \alpha^a (J)$.
Hence $\alpha^a: \setmono\to \setmono$ is a squarefree stable operator.
\end{proof}

\begin{cor}\label{exball}
Let $n>d>0$ be positive integers, $B_d(U)$ a squeezed $d$-ball on $[n]$ and
$\alpha(I(U))\subset \aru n$ the ideal generated by all monomials $\alpha(u)$ with $u\in\aru {n-d-1}$ and with $u\not\in U$.
Then
\begin{itemize}
\item[(i)]
$I_{\dele {B_d(U)}} = \alpha(I(U))$.
\item[(ii)]
$\dele {B_d(U)}$ is the simplicial complex generated by
$$L=\{ \{i_1,i_2+1,\dots,i_k+k-1\}\cup \{n-k+1,n-k+2,\dots,n\}:x_{i_1}x_{i_2}\cdots x_{i_k}\in U\}.$$ 
\end{itemize}
\end{cor}

\begin{proof}
(i)
Let $I(U) \subset \aru n$ be the ideal generated by all monomials $u \in \aru {n-d-1}$ with $u \not\in U$
and $\alpha^2:\setmono \to \setmono$ the stable operator defined in \S 4.
Then Proposition \ref{ginball} says that
$I_{B_d(U)}$ is the ideal generated by all monomials $\alpha^2(u)=\alpha (\alpha (u))$ with $u \in \aru {n-d-1}$ and with 
$u\not \in U$.
Since Lemma \ref{sss} says that $\alpha(I(U))$ is squarefree strongly stable
and since  $\alpha: \setmono \to \setmono$ is a squarefree stable operator,
it follows from Proposition \ref{sqoperator} that
$$I_{\dele {B_d(U)}}= \alpha(I(U)).$$

(ii)
First, we will show $L \subset \dele {B_d(U)}$.
Let  $u=x_{i_1}x_{i_2}\cdots x_{i_k}\in U$ with $i_1 \leq i_2 \leq \cdots \leq i_k$.
The ideal $\alpha(I(U))$ is squarefree strongly stable.
Then we have $\alpha(u) \not\in \alpha (I(U))$,
because if $ \alpha(u) \in \alpha (I(U))$ then there is
$ \alpha(w) = x_{i_1} x_{i_2+1}\cdots x_{i_l+{l-1}} \in G(\alpha (I(U)))$
for some $l \leq k$ and $w \not \in U$ divides $u$.
Since $I_{\dele {B_d(U)}}=\alpha (I(U))$, we have
$\{i_1 , i_2+1, \dots,i_{k}+k-1\}  \in  \dele {B_d(U)}$.
Since $B_d(U)$ is Cohen--Macaulay,
it follows from \cite[Theorem 8.13]{H} that $\dele {B_d(U)}$ is pure.
Thus there is a $(d-k)$-subset $F \subset [n] \setminus \{i_1 , i_2+1, \dots,i_{k}+k-1\} $ such that
$F \cup \{i_1 , i_2+1, \dots,i_{k}+k-1\} \in B_d(U)$.
Since $\dele {B_d(U)}$ is shifted, we may assume $F= \{n-k+1,n-k+2,\dots,n\}$.
Thus we have $L \subset \dele {B_d(U)}$.

On the other hand,
Lemma \ref{fosb} says
$$f_d(\dele {B_d(U)}) = f_d(B_d(U)) =\sum_{j=0}^{d+1} h_j(B_d(U))=|U|=|L|.$$
Since $\dele {B_d(U)}$ is pure and since $L \subset \dele {B_d(U)}$,
it follows that $\dele {B_d(U)}$ is the simplicial complex generated by $L$.
\end{proof}

\noindent
\textbf{Acknowledgement}:
%
I would like to thank the referee for correcting the proof of Lemma \ref{4polytopal} of the original manuscript,
encouraging me to find Example \ref{referee},
and many helpful suggestions to improve the paper.

\end{document}